\documentclass[12pt]{amsart}
\usepackage{amsmath}
\usepackage{amsfonts}
\usepackage{amssymb}
\usepackage{epsfig,latexsym,graphicx,stmaryrd}

\newtheorem{Thm}{Theorem}[section]
\newtheorem{thm}[Thm]{Theorem}

\newtheorem{corollary}[Thm]{Corollary}
\newtheorem{lemma}[Thm]{Lemma}

\newtheorem{remark}[Thm]{Remark}

\newtheorem{alg}[Thm]{Algorithm}

\newcommand{\PP}{\mathbb P}

\newcommand{\TTT}{{\mathcal{T}}}

\newcommand{\HH}{{\cal H}}

\newcommand{\ip}[2]{{\langle#1,#2\rangle}}
\newcommand{\bigip}[2]{{\left\langle#1,#2\right\rangle}}
\newcommand{\norm}[1]{{\|#1\|}}
\newcommand{\outp}[2]{{\llbracket #1,#2 \rrbracket }}

\def\ldots{\mathinner{\ldotp\ldotp\ldotp}}
\def\ldots{\mathinner{\cdotp\cdotp\cdotp}}

\def \K{\rm {\bf K}}

\def \cal{\mathbb}

\def \beq{\begin{eqnarray*}}
\def \eeq{\end{eqnarray*}}

\def \R{\mathbb{R}}
\def \C{\mathbb{C}}
\def \K{\mathbb{K}}
\def \E{\mathbb{\cal E}}

\newcommand{\fc}{{\mathbb F}}

\newcommand{\Ac}{\mathcal{ A}}

\newcommand{\nonlin}{{\alpha}}
\renewcommand{\SS}{\mathcal{S}}
\newcommand{\Sonezero}{{\SS^{1,0}}}
\newcommand{\Soneone}{{\SS^{1,1}}}
\newcommand{\Phit}{{\Phi'}}

\newcommand{\HHH}{K}

\newcommand{\epsi}{\epsilon}
\newcommand{\Gcal}{{\Upsilon}}

\renewcommand{\i}{{\bf j}}

\title[Reconstruction from Magnitudes of Frame Coefficients]{Reconstruction of Signals from Magnitudes of Redundant Representations: The Complex Case}
\author[R.~Balan]{Radu Balan}
\address{Department of Mathematics \\ University of Maryland, College Park MD 20742}
\email[R.~Balan]{rvbalan@math.umd.edu} 

\begin{document}

\begin{abstract}
This paper is concerned with the question of reconstructing
a vector in a finite-dimensional complex Hilbert space when 
only the magnitudes of the coefficients of the vector under a redundant linear map
are known. We present new invertibility results as well as an iterative algorithm that finds the least-square solution
which is robust in the presence of noise. We analyze its numerical performance by
comparing it to the Cramer-Rao lower bound.
\end{abstract}

\maketitle

{\bf AMS (MOS) Subject Classification Numbers}: 15A29, 65H10, 90C26
\\
\\

{\bf Key words:} Frame, phase retrieval, Cramer-Rao lower bound, phaseless reconstruction
\\
\\

\section{Introduction}

This paper is concerned with the question of reconstructing
a vector $x$ in a finite-dimensional {\it complex } Hilbert space $H$ when 
only the magnitudes of the coefficients of the vector under a redundant linear map
are known.

Specifically our problem is to reconstruct $x\in H$ up to a global phase factor from the magnitudes 
$\{ |\ip{x}{f_k}|~,~1\leq k\leq m\}$ where $\{f_1,\ldots,f_m\}$ is a frame (complete system) for $H$.
The real case was considered in \cite{Bal12a}. Here we develop the theory for the complex case.
 
A previous paper \cite{BCE06} described the importance of this problem
to signal processing, in particular to the analysis of speech. The problem appears in X-Ray crystallography 
under the name of "phase retrieval" problem (see \cite{reviewBates}) where the frame vectors are the Fourier frame vectors.
The case of the windowed Fourier transform was considered in the '80s \cite{NQL82}; see also \cite{Bal2010} for a frame based approach to this case.
A different approach is taken in \cite{Raz13}. The authors propose a novel algorithm
adapted to compactly supported signals and FFT that uses individual signal spectral
powers and two additional interferences between signals. A 3-term polarization identity has been used in \cite{ABFM12}
together with the angular synchronization algorithm. 

%A different approach (which we called the ''algebraic approach'') was proposed
%in \cite{Bal2009}. While it applies to both real and complex cases, noisless
%and noisy cases, the approach requires solving a linear system of exponential size 
%in space dimension.  The algebraic approach mentioned earlier
%generalizes the approach in \cite{BBCE07} where reconstruction is performed
%with complexity $O(n^2)$ (plus computation of the principal eigenvector for a
%matrix of size $n$) but requiring $m=O(n^2)$ frame vectors.

Recently the authors of \cite{CSV12} developed a convex
optimization algorithm (PhaseLift) and proved its ability to perform exact reconstruction 
in the absence of noise, as well as its stablity under noise conditions. 
In a separate paper, \cite{CESV12}, the authors developed further a similar algorithm in the case
of windowed DFT. The original requirement of $m=O(n\,log(n))$ vectors has been
relaxed to $m=O(n)$ in \cite{CL12}. Similar convex optimization solutions have been proposed by the authors
of \cite{DemHan12} and \cite{Wald12}. Additionally, \cite{Hand13} studied the duality gap in this approach
and obtained a necessary and sufficient condition for the existence of a dual certificate.

While writing this paper, we become aware of \cite{BCMN13} where the authors obtained similar results to the 
injectivity criteria presented here, as well as to the Cramer-Rao Lower Bound derived in this paper. We will comment more later in the paper.
We also acknowledge the paper \cite{EldarMend12} where certain Lipschitz bounds have been obtained in the real case.
The real case of the stabilty bounds obtained here are presented in a separate paper \cite{BaWa13a} together with 
additional results for the real case.

The organization of the paper is as follows. In section \ref{sec2} we define the problem
explicitly. In section \ref{sec3} we describe new analysis results. Specifically we analyze in more detail spaces of symmetric
operators of constrained signature, and then we show how they are related to the phaseless recovery problem.
Our results are canonical, meaning they are independent to a particular choice of basis. 
In section \ref{sec4} we establish two robustness results: bi-Lipschitzianity of the nonlinear analysis map, and the Cramer-Rao Lower Bound (CRLB).
In section \ref{sec5} we present a new reconstruction algorithm based on the Least-Square method. We also obtain robustness bounds to noise.
Its performance is analyzed in section \ref{sec6} and is compared to the CRLB. Section \ref{sec7} contains
conclusions and is followed by references.

\section{Background\label{sec2}}

Let $H$ be an n-dimensional complex Hilbert space (such as $\C^n$ or $\C^{p_1\times p_2}$ the vector space of $p_1\times p_2$ complex matrices)  with scalar product $\ip{}{}$ linear in the first term
and antilinear in the second term and a {\em conjugation} $c:H\rightarrow H$ (see e.g. \cite{Halmos74}; conjugation is an antilinear transformation that squares to the identity). 
Let $\fc=\{f_1,\ldots,f_m\}$ be a spanning set of $m$ vectors in $H$. As $H$ has finite dimension
such a set forms a {\em frame}. In the infinite
dimensional case, the concept of frame involves a stronger property than completeness (see 
for instance \cite{Cass-artofframe}). 
We review additional terminology and properties which remain still true in the infinite dimensional setting.
The set $\fc$ is a frame if and only if there are two positive constants $0<A\leq B<\infty$ (called frame bounds)
such that
\[ A\norm{x}^2 \leq \sum_{k=1}^m |\ip{x}{f_k}|^2 \leq B\norm{x}^2 ~~.\]
When we can choose $A=B$ the frame is said to be {\em tight}. For $A=B=1$ the frame is called {\em Parseval}.
A set of vectors $\fc$ of the $n$-dimensional Hilbert space $H$ is said to have {\em full spark} if every subset
of $n$ vectors is linearly independent (see \cite{ACM12} for a full discussion of such frames). 

\subsection{Problem Definition and Notations}

For a vector $x\in H$, the collection of coefficients
$\{ \ip{x}{f_k}~,~1\leq k\leq m\}$ represents the analysis of the vector $x$ given by the frame $\fc$.
In $H$ we consider the following equivalence relation:
\begin{equation}
\label{eq:equiv}
x,y\in H ~~,~~ x\sim y~~{\rm iff}~~y=zx~{\rm for~some~scalar~z~with}~|z|=1.
\end{equation}
Note $z=e^{i\varphi}$ for some real number $\varphi$.
Let $\hat{H}=H/\sim$ be the set of classes of equivalence induced by this relation.
Thus $\hat{x}=\{e^{i\alpha}x,
0\leq\alpha <2\pi \}$. The analysis map induces
the following nonlinear map
\begin{equation}
\label{eq:varphi}
\nonlin :\hat{H}\rightarrow (\R^{+})^m~~,~~\nonlin(\hat{x})=(|\ip{x}{f_k}|^2)_{1\leq k\leq m}
\end{equation}
where $\R^{+}=\{x~~,~~x\in\R~,~x\geq 0\}$ is the set of nonnegative real numbers.
In \cite{Bal12a} we studied when the nonlinear map $\nonlin$ is injective, mostly in the real case, and we provided some
necessary conditions of injectivity in the complex case. We review these results below.
In this paper we obtain additional injectivity results. We then concentrate on the additive white Gaussian noise model
\begin{equation}
\label{eq:model}
y = \nonlin(x) + \nu ~~,~~\nu\sim \mathcal{N}(0,\sigma^2).
\end{equation}
We describe an algorithm (the Iterative Regularized Least-Square (IRLS) algorithm) to solve the estimation problem. 
We prove some convergence results and we study its performance in the noisy case.
We shall derive the Cramer-Rao Lower Bound (CRLB) for this model and compare the algorithm performance to this bound.

We describe several objects that will be used in this paper. 

The set $B(H)$ denotes the set of bounded linear operators on $H$. In $B(H)$ for any $1\leq p\leq\infty$, the $p$-norm
of $T\in B(H)$ denoted $\norm{T}_p$ is given by the $p$-norm of its set of singular eigenvalues. In particular 
for $p=1$, the 1-norm $\norm{T}_1$ is called the nuclear norm of $T$; for $p=2$, the 2-norm $\norm{T}_2=\sqrt{tr(T^*T)}$
is the Frobenius norm of $T$; for $p=\infty$, the $\infty$-norm $\norm{T}_{\infty}$ is the same as the operator norm of $T$ on $H$,
simply denoted $\norm{T}$.

In general, for two operators $T,S\in B(H)$ we denote $\ip{T}{S}_{B(H)} = tr\{TS^*\}$ their
Hilbert-Schmidt scalar product, where $S^*$ is the adjoint of $S$ and $tr$ denotes the trace. Note the scalar product is basis independent, but it depends on the underlying Hilbert space structure.
When no danger of confusion we drop the index $B(H)$ from the scalar product notation.

For each frame vector $f_k$ we denote by $F_k$ its associated rank-1 operator
\begin{equation}
\label{eq:Fk}
F_k: H\rightarrow H~~,~~F_k(x) = \ip{x}{f_k}f_k.
\end{equation}
In general the {\em associated rank-1 operator} to a vector $x\in H$ is the operator $X:H\rightarrow H$, $X=xx^*$ which acts by $X(v)=\ip{v}{x}x$.
Here and throughout the paper $x^*$ denotes the adjoint (or dual) of $x$, that is $x^*:H\rightarrow \C$, $x^*(v)=\ip{v}{x}$.
Note 
%this is a slight misnomer: 
$X$ has {\em at most} rank one. Specifically, $X$ has rank one if and only if $x\neq 0$;
otherwise $X$ has rank zero.
%One can think of $X$ as simply the $n\times n$ matrix $xx^*$, where $x^*$ denotes the transpose and conjugate of $x$.
%However, as much as possible, we plan to use canonical notations (that is, basis independent).
 
For any two vectors $u,v\in H$ we define their {\em symmetric outer product} denoted $\outp{u}{v}$ by
\begin{equation}
\label{eq:outp}
\outp{u}{v}:H\rightarrow H ~,\outp{u}{v}=\frac{1}{2}(uv^*+vu^*)~,~\outp{u}{v}(x) =\frac{1}{2} (\ip{x}{u}v + \ip{x}{v}u)~~.
\end{equation}
%One can think of $\outp{u}{v}$ as the $n\times n$ symmetric matrix $\frac{1}{2}(vu^*+uv^*)$.
Note the rank-1 operator associated to a vector $x$ can be written as $\outp{x}{x}$. In particular $F_k=\outp{f_k}{f_k}$. 
Note also $\outp{u}{v}$ is $\R$-bilinear
but it is not $\C$-(bi)linear. Furthermore $\outp{u}{v}=\outp{v}{u}$.

Following \cite{BBCE07} the nonlinear map $\nonlin$ induces a linear map
$\Ac$ on the set $B(H)$ of bounded operators on $H$:
\begin{equation}
\label{eq:Ac}
\Ac: B(H) \rightarrow \C^m ~~,~~(\Ac(T))_k = \ip{Tf_k}{f_k}=tr\{TF_k\}~,~1\leq k\leq m.
\end{equation}
Thus $(\Ac(T))_k=\ip{T}{F_k}_{B(H)}$. Also note $\nonlin(x)=\Ac(X)$ where $X=xx^*$
is the rank-1 operator associated to $x$. This remark was first observed by B. Bodmann in \cite{BBCE07}.

Let $Sym(H)$ denote the set of self-adjoint operators on $H$, $Sym(H)=\{T\in B(H)~,~T^*=T\}$. We denote by $\SS^{p,q}$ or $\SS^{p,q}(H)$,  the set of
self-adjoint operators on $H$ that have at most $p$ positive eigenvalues and at most $q$ negative eigenvalues:
\begin{equation}
\label{eq:Spq}
\begin{split}
\SS^{p,q} = &  \{ T\in Sym(H)~,~Sp(T)=\{\lambda_1,\ldots,\lambda_n\}~,~ \\
 & \lambda_1\geq\cdots\geq \lambda_p \geq 0=\lambda_{p+1} =\cdots
=\lambda_{n-q} \geq \lambda_{n-q+1} \geq \cdots\geq \lambda_n \}
\end{split}
\end{equation}
where $Sp(T)$ denotes the spectrum of $T$ (i.e. the set of its eigenvalues).
Notice $\SS^{p,q}$ is not a linear space, but instead it is a cone in $B(H)$. This cone property is key in deriving 
robustness and stability bounds later on.

 We denote by $\lambda_{max}(T)$ the largest eigenvalue of $T$ and by $\lambda_{min}(T)$
the smallest eigenvalue of $T$. 
 In particular we are interested in $\Sonezero$ and $\Soneone$:
\begin{eqnarray}
\label{eq:S10}
\Sonezero & = & \{T\in Sym(H)~,~rank(T)\leq 1,\lambda_{min}(T)=0 \} \\
\label{eq:S11}
\Soneone & = & \{ T\in Sym(H)~,~rank(T)\leq 2,Sp(T)=\{\lambda_{max}(T),0^{(n-2)},\lambda_{min}(T) \}, \\
 & & \lambda_{max}(T)\geq 0\geq \lambda_{min}(T) \} \nonumber
\end{eqnarray}
Note the following obvious inclusions
\begin{equation}
\{ 0\} \subset \Sonezero\subset\Soneone \subset Sym(H)   ~,~\{ 0\} \subset \SS^{0,1}\subset\Soneone \subset Sym(H)
\end{equation}
We denote by $\mathring{\SS}^{p,q}$ the subset of $\SS^{p,q}$ of selfadjoint operators that have rank $p+q$, hence
exactly $p$ strictly positive eigenvalues and $q$ strictly negative eigenvalues. Thus
\begin{equation}
\label{eq:Spqdecomp}
\SS^{p,q} = \mathring{\SS}^{p,q} \cup \SS^{p-1,q} \cup \SS^{p,q-1}
\end{equation}
represents a disjoint partition of $\SS^{p,q}$. In particular
\begin{equation}
\Soneone = \mathring{\SS}^{1,1} \cup \SS^{0,1} \cup \SS^{1,0} ~,~\Sonezero = \mathring{\SS}^{1,0}\cup \{0\}~.
\end{equation}
Finally we let $GL(H)$ denote the group of invertible linear operators on $H$.  
A more detailed analysis of these sets is presented in subsection \ref{subsec3.2}.

Next we describe the realification of the Hilbert space $H$. To do so canonically we need to fix a conjugation $c:H\rightarrow H$. 
To the complex Hilbert space $H$ with conjugation $c$ we associate its $2n$-dimensional  real vector space $H_\R$ subset of $H\times H$ built from vectors
$v_R = \frac{1}{2}(v+c(v))$ and $v_I = \frac{1}{2i}(v-c(v))$ as follows:
\begin{equation}
\label{eq:Hreal}
H_{\R} = \left\{ \left( \frac{1}{2}(v+c(v)) ,\frac{1}{2i} (v-c(v)) \right)~,~v\in H \right\}.
\end{equation}
Thus $H_\R$ is the image of $H$ through the $\R$-linear map,
\begin{equation}
\i:H\rightarrow H\times H~,~\i(v) = \left( \frac{1}{2}(v+c(v)) ,\frac{1}{2i} (v-c(v)) \right).
\end{equation}
Note $\i$ is injective with range $H_\R$. Furthermore $\i:H\rightarrow H_{\R}$ is a norm preserving $\R$- isomorphism. 
Its inverse is given by
\begin{equation}
\i^{-1}:H_\R \rightarrow H~~,~~\i^{-1}(u,v) = u + iv~.
\end{equation}
Let $J$ denote the linear map defined by 
\begin{equation}\label{eq:J}
J:H\times H\rightarrow H\times H~,~J(v,w) = (-w,v)~.
\end{equation}
Note it is conjugate to the multiplication by $i$ in $H$:
\begin{equation}
J:H_\R\rightarrow H_\R ~,~J(\i(v)) = \i(iv)~.
\end{equation}
Hence $H_\R$ is $J$ invariant.
In $H_\R$ the induced scalar product is given by
\begin{equation}
\ip{\i(v)}{\i(w)} := \bigip{\frac{1}{2}(v+c(v))}{\frac{1}{2}(w+c(w))} + \bigip{\frac{1}{2i}(v-ic(v))}{\frac{1}{2i}(w-c(w))} = real(\ip{v}{w}).
\end{equation}
We denote by $\ip{v}{w}_\R = real(\ip{v}{w})$ the $\R$-linear inner product on $H$. Note
\begin{equation}
\label{eq:ip}
\begin{split}
 \ip{u}{v} & =real(\ip{u}{v}) + i\,imag(\ip{u}{v}) = \ip{u}{v}_{\R} - i \ip{iu}{v}_{\R} \\ 
 & = \ip{u}{v}_{\R} + i \ip{u}{iv}_{\R} =\ip{\i(u)}{\i(v)} + i \ip{\i(u)}{J\,\i(v)}  .
\end{split}
\end{equation}
If $\{e_1,\ldots,e_n\}$ is an orthonormal basis in $H$, then $\{\i(e_1),\ldots,\i(e_n),J\i(e_1),\ldots,J\i(e_n)\}$ is an orthonormal basis in $H_\R$.
Which shows the real dimension of $H_\R$ is $2n$, $\dim_{\R}H_\R = 2n$. 

On $H\times H$ there are two inner product structures:
\begin{eqnarray}
\ip{(x,y)}{(u,v)}_{\C} & = & \ip{x}{u}+\ip{y}{v} \\
\ip{(x,y)}{(u,v)}_{\R} & = & \frac{1}{2}\left( \ip{(x,y)}{(u,v)} + \ip{(u,v)}{(x,y)} \right) = \ip{x}{u}_{\R} + \ip{y}{v}_\R
\end{eqnarray}
With respect to $\ip{\cdot}{\cdot}_{\C}$, $H\times H$ is a $\C$-vector space of dimension $2n$. With respect to $\ip{\cdot}{\cdot}_{\R}$, $H\times H$
is a $\R$-vector space of dimension $4n$. On $H_\R$ the two inner products coincide. Because of this fact we shall simply denote $\ip{\xi}{\eta}_{H_\R}$ or $\ip{\xi}{\eta}$ the scalar product
on $H\times H$ whenever $\xi,\eta\in H_\R$. Furthermore $\i$ is norm preserving $\norm{x}_H=\norm{\i(x)}_{H_\R}$ and 
$\ip{x}{y}_\R = \ip{\i(x)}{\i(y)}_{H_{\R}}$ for all $x,y\in H$.
The orthogonal complement of $H_\R$ in $H\times H$ with
respect to $\ip{\cdot}{\cdot}_{\R}$ is given by
\begin{equation}
\label{eq:HRperp}
H_\R^{\perp} = iH_\R = \left\{ \left( \frac{i}{2}(x+c(x)),\frac{1}{2}(x-c(x)) \right)~,~x\in H  \right\}.
\end{equation}
The orthogonal projection onto $H_\R$ with respect to the real structure is given by
\begin{equation}
P_\R: H\times H\rightarrow H_\R\subset H\times H ~~,~~P_\R(u,v) = \left( \frac{1}{2}(u+c(u)),\frac{1}{2}(v+c(v)) \right).
\end{equation}
The projection onto the orthogonal complement $H_\R^{\perp}$ is given by
\begin{equation}
 P_\R^{\perp} = 1-P_\R~~,~~P_\R^{\perp}(u,v)=\left(\frac{1}{2}(u-c(u)),\frac{1}{2}(v-c(v)) \right) = i \left(
\frac{1}{2i}(u-c(u)),\frac{1}{2i}(v-c(v))\right).
\end{equation}
Note:
\begin{equation}
P_\R^{\perp}(u,v) = -i P_\R(iu,iv)~.
\end{equation}
Fix $a,b\in H$ and define the linear operator on $H$,
\begin{equation}
T_{a,b}:H\rightarrow H~~,~~T_{a,b}(x) = \ip{x}{a}b~.
\end{equation}
Associate the following $\R$-linear operator on $H\times H$
\begin{equation}
\label{eq:Ttildeab}
\tilde{T}_{a,b}:H\times H\rightarrow H_\R\subset H\times H ~~,~~ \tilde{T}_{a,b}(u)=\ip{u}{\i(a)}_{\R} \i(b) + \ip{u}{\i(ia)}_{\R} \i(ib)~.
\end{equation}
Note $\ip{\i(b)}{\i(ib)}_{\R}=0$ and $H_\R$ is invariant under the action of $T_{a,b}$.
Direct computations show the following diagram is commutative:
\begin{equation}
\begin{array}{cccccc}
\mbox{$T_{a,b}$} & \mbox{$:$} & \mbox{$H$} & \mbox{$\stackrel{T_{a,b}}{\longrightarrow}$} & \mbox{$H$} & \\
& \mbox{$$} & \mbox{$\i\downarrow$} & & \mbox{$\downarrow\i$} & \mbox{$$} \\
\mbox{$\tilde{T}_{a,b}$} & \mbox{$:$} & \mbox{$H_\R$} & \mbox{$\stackrel{\tilde{T}_{a,b}}{\longrightarrow}$} & \mbox{$H_\R$} & 
\end{array}
\end{equation}
Similarly each symmetric operator $T$ in $Sym(H)$ gets mapped into  a symmetric operator in $Sym(H_\R)$ of double rank. Denote by $\tau$ this mapping, 
$\tau:Sym(H)\rightarrow Sym(H_\R)$ that makes the following diagram commutative:
\begin{equation}
\label{eq:tauT}
\begin{array}{ccc}
\mbox{$H$} & \mbox{$\stackrel{T}{\longrightarrow}$} & \mbox{$H$} \\
\mbox{$\i\downarrow$} & & \mbox{$\downarrow\i$} \\
\mbox{$H_\R$} & \mbox{$\stackrel{\tau(T)}{\longrightarrow}$} & \mbox{$H_\R$}  
\end{array}
\end{equation}
If desired, $\tau(T)$ can be extended to $H\times H$ using the $\R$-linear scalar product $\ip{\cdot}{\cdot}_{\R}$ on $H\times H$. 
$\tau$ is $\R$-linear but  not $\C$-linear. 
In particular:
\begin{eqnarray}
T=\outp{x}{x}\in\Sonezero & \Longmapsto & \tau(T)=\outp{\i(x)}{\i(x)}+\outp{\i(ix)}{\i(ix)}\in\SS^{2,0}(H_\R) \label{eq:P}\\
T= \outp{x}{y}\in\Soneone & \Longmapsto & \tau(T)=\outp{\i(x)}{\i(y)}+\outp{\i(ix)}{\i(iy)}\in\SS^{2,2}(H_\R)  \label{eq:PQ}
\end{eqnarray}
Using the $\R$-linear operator $J$ introduced in (\ref{eq:J}) the first relation above can be rewritten as
\begin{equation}
\tau(xx^*) = \xi\xi^* + J\xi\xi^*J^*~~,~~{\rm where}~\xi=\i(x)\in H_\R
\end{equation}
and the adjoints $\xi^*$ and $J^*$ are taken with respect to the scalar product of $H_\R$.
In general an operator in $\SS^{p,q}(H)$ gets mapped into an operator in $\SS^{2p,2q}(H_\R)$.
Note the map $\tau$ preserves scalar products between selfadjoint operators up to a multiplicative constant: 
\begin{equation}
\ip{T}{S}_{B(H)} = 2\ip{\tau(T)}{\tau(S)}_{B(H_\R)}~.
\end{equation}
Thus $\tau$ is a monomorphism (i.e. linear injective morphism) between the $\R$-vector spaces $Sym(H)$ and $Sym(H_\R)$. 
We shall discuss in more details the linear map $\tau$ in subsection \ref{subsec3.21}.
\vspace{3mm}

Note, in the case $H=\C^n$, $c(z)=(\bar{z}_1,\ldots,\bar{z}_n)^T$, and $\ip{v}{w}=v^Tc(w)$,
the map $\i$ acts by $z\in\C^n\mapsto \i(z) = (real(z^T),imag(z^T))^T$.  Here $\mbox{}^T$ denotes transposition.
Thus $H_\R=\{(real(v^T),imag(v^T))^T~,~v\in\C^n\}=\R^{2n}$, and the $\R$-linear map $J$ has the block form
\[ J=\left[ \begin{array}{cc} 
\mbox{$0$} & \mbox{$-I$} \\
\mbox{$I$} & \mbox{$0$}
\end{array} \right]
\]
where $I$ is the identity matrix of size $n$. The scalar product in $H_\R$ is the usual real scalar product. For any vector $\xi\in H_\R=\R^{2n}$  
the adjoint reduces to transposition $\xi^*=\xi^T$. 
\vspace{3mm}

We return to the frame set $\fc=\{f_1,\ldots,f_m\}$ and the Hilbert space $H$.
We let 
\begin{equation}
\Phi_k = \tau(f_k f_k^*) = \varphi_k \varphi_k^* + J\varphi_k\varphi_k^*J^*
\end{equation}
denote the image of $F_k=\outp{f_k}{f_k}$ under $\tau$, where $\varphi_k=\i(f_k)$ is an element of $H_\R$.

The last notation we introduce here is the following map on $H_\R$:
\begin{equation}
\label{eq:I}
R : H_\R \rightarrow Sym(H_\R) ~~,~~R(\xi) = \sum_{k=1}^m \Phi_k \xi\xi^* \Phi_k.
\end{equation}
As we will see  later $R(\xi)$ is related to the Fisher information matrix for the measurement model (\ref{eq:model}), and it plays a key role 
in obtaining a necessary and sufficient condition of injectivity for the nonlinear map $\nonlin$. More explicit
%\[ 
%\Phi_k \outp{\xi}{\xi} \Phi_k = \outp{\Phi_k \xi}{\Phi_k \xi}
%\]
%$R(\xi)$ has the following more explicit form:
\begin{equation}
\label{eq:Iexp}
R(\xi) = \sum_{k=1}^m v_kv_k^*~~, v_k=\Phi_k\xi=\ip{\xi}{\varphi_k}\varphi_k + \ip{\xi}{J\varphi_k}J\varphi_k.
\end{equation}
We shall not overload the notation and use the same letter $R$ to denote the map $R:H\rightarrow Sym(H_\R)$, defined by $x\mapsto R(\i(x))$.
Finally, we let $\delta_{i,j}$ denote the Kroneker symbol: $\delta_{i,j}=1$ if $i=j$, and $0$ otherwise.

\subsection{Existing Results}

We revise now existing results on injectivity of the nonlinear map $\nonlin$. A subset $Z$ of a topological space $(X,{\cal \tau})$ 
is said to be {\em generic} if its open interior is dense. In the following statements, the term
{\em generic} refers to the Zarisky topology: a set $Z\subset \K^{n\times m}= \K^{n}\times\cdots\times \K^n$ is said to be {\em generic}
if $Z$ is dense in $\K^{n\times m}$ and its complement is a finite union of zero sets of polynomials in $nm$ variables 
with coefficients in the field $\K$ (here $\K=\R$ or $\K=\C$). 
\begin{thm} \label{th2.1} In the real case when $H=\R^n$ the following are equivalent:
\begin{enumerate}
\item The nonlinear map $\nonlin$ is injective;
\item (\cite{BCE06}, Th.2.8) For any disjoint partition of the frame set $\fc=\fc_1\cup\fc_2$, either $\fc_1$ spans $H$ or $\fc_2$ spans $H$.

\item (\cite{Bal12a},Th.2.4(2)) For any two vectors $x,y\in H$ if $x\neq 0$ and $y\neq 0$ then 
$$\sum_{k=1}^m |\ip{x}{f_k}|^2|\ip{y}{f_k}|^2 >0$$

\item (\cite{Bal12a},Th.2.4(3)) There is a real constant $a_0>0$ so that for all $x,y\in H$,
\begin{equation}
\label{eq:q4}
\sum_{k=1}^m |\ip{x}{f_k}|^2 |\ip{y}{f_k}|^2 \geq a_0 \norm{x}^2 \norm{y}^2
\end{equation}

\item (\cite{Bal12a},Th.2.4(4)) There is a real constant $a_0>0$ so that for all $x\in H$,
\begin{equation}
\label{eq:q2}
R(x) := \sum_{k=1}^m |\ip{x}{f_k}|^2 \ip{\cdot}{f_k}f_k \geq a_0 I
\end{equation}
where the inequality is in the sense of quadratic forms.
\end{enumerate}

Additionally, the following statements hold true:
\begin{enumerate}
\item (\cite{BCE06}, Prop2.5) If $\nonlin$ is injective then  $m\geq 2n-1$;

\item (\cite{BCE06}, Prop.2.5) If $m\leq 2n-2$ then $\nonlin$ cannot be injective;

\item (\cite{BCE06}, Cor.2.7(1)) If $m=2n-1$ then $\nonlin$ is injective if and only if $\fc$ is full spark;

\item (\cite{BCE06}, Cor.2.6) If $m\geq 2n-1$ and $\fc$ is full spark then the map $\nonlin$ is injective;

\item (\cite{BCE06}, Th.2.2)  If $m\geq 2n-1$ then for a generic frame $\fc$ the map $\nonlin$ is injective.
\end{enumerate}
\end{thm}

In the complex case the following results are known:

\begin{thm} \label{th2.3} In the complex case when $H=\C^n$ the following statements hold true:
\begin{enumerate}
\item (\cite{BCE06}, Th.3.3) If $m\geq 4n-2$ then for a generic frame $\fc$ the map $\nonlin$ is injective;
%\item (\cite{Fink04}Th.II) If $\nonlin$ is injective then $m\geq 3n-2$;
%\item (\cite{Fink04}Th.II) If $m\leq 3n-3$ then the map $\nonlin$ cannot be injective.
\item (\cite{BH13}) For any positive integer $n$ there is a frame with $m=4n-4$ elements so that the nonlinear map $\nonlin$ is injective;
\item (\cite{HMW11}, Corollary 4) If $\nonlin$ is injective then 
\begin{equation}
\label{eq:HMW11}
m\geq 4n-2-2\beta + \left\{ 
\begin{array}{rl}
\mbox{$2$} & \mbox{if $n$ odd and $\beta=3\,mod\,4$} \\
\mbox{$1$} & \mbox{if $n$ odd and $\beta=2\,mod\,4$} \\
\mbox{$0$} & \mbox{otherwise}
\end{array} \right.
\end{equation}
where $\beta=\beta(n)$ denotes the number of 1's in the binary expansion of $n-1$.
\item The following are equivalent:
\begin{enumerate}
\item The map $\nonlin$ is injective;
\item (\cite{HMW11} Prop. 2) 
\begin{equation}
\label{eq:kerAc}
ker(\Ac) \cap \left(\Sonezero - \Sonezero \right) = \{ 0\}
\end{equation}
\item(\cite{BCMN13}, Theorem 4) $dim\, S(u)\geq 2n-1$ for every $u\in\C^n\setminus\{0\}$;
\item (\cite{BCMN13}, Theorem 4) $S(u)=span_{\R}(iu)^{\perp}$ for every $u\in\C^n\setminus\{0\}$.
\item (\cite{CEHV15}\footnote{This result was not known at the time the present paper was submitted for publication.}, Theorem 1.1) If $m\geq 4n-4$ then for a generic frame $\fc$ the map $\alpha$ is injective;
\end{enumerate}
\end{enumerate}
\end{thm}

In the last two conditions, $S(u)=span_{\R} \{f_k f_k^* u\}_{1\leq k\leq m}$, where $f_kf_k^* u$ is seen as a $2n$ vector in $\R^{2n}$.

\section{New Analysis Results\label{sec3}}
This section contains our injectivity results of the nonlinear map $\nonlin$ as well as an in-depth analysis of the spaces $\SS^{p,q}$.

\begin{thm}
\label{t3.1}
Let $H$ be a $\C$-vector space of dimension $n$, with scalar product $\ip{}{}$ and conjugation $c:H\rightarrow H$. The following are equivalent:
\begin{enumerate}
\item The nonlinear map $\nonlin:\hat{H}\rightarrow\R^m$, $(\nonlin(x))_k=|\ip{x}{f_k}|^2$ is injective.
\item There is a constant $a_0>0$ so that for every $u,v\in H$ 
\begin{equation}
\label{eq:l2bound}
\sum_{k=1}^m |\ip{F_k}{\outp{u}{v}}|^2 \geq a_0 \norm{\outp{u}{v}}_1^2
\end{equation}
where $F_k=f_kf_k^*$. Explicitly, this means:
\begin{equation}
\label{eq:l2bound2}
\sum_{k=1}^m \left( real(\ip{u}{f_k}\ip{f_k}{v}) \right)^2 \geq a_0 \left[ \norm{u}^2 \norm{v}^2 - (imag(\ip{u}{v}))^2\right] .
\end{equation}
\item For any $\xi\in H_\R$, $\xi\neq 0$, $rank(R(\xi)) = 2n-1$.
\item There is $a_0>0$ so that for all $\xi\in H_\R$, $\xi\neq 0$,
\begin{equation}
\label{eq:Ibound}
R(\xi) \geq a_0 \norm{\xi}^2 P_{J\xi}^{\perp}
\end{equation}
where the inequality holds in the sense of quadratic forms in $H_\R$, and
\begin{equation}
P_{J\xi}^{\perp} = 1 - \frac{1}{\norm{\xi}^2 }J\xi\xi^*J^*
\end{equation}
is the orthogonal projection in $H_\R$ onto the orthogonal complement of $J\xi$ in $H_\R$.
\end{enumerate}
\end{thm}

The proof is given in section \ref{subsec3.3}.

\begin{remark}
The two constants $a_0$ in (2) and (4) can be chosen to be equal, hence the same notation. We will see in the next section,
this common constant is related to robustness and stability of any reconstruction scheme.
\end{remark}

\begin{remark}
The proof of $(3)\Leftrightarrow(4)$ shows that the optimal bound $a_0$ is given by
\begin{equation}\label{eq:a0}
a_0^{opt} = min_{\xi\in H_\R,\norm{\xi}=1} a_{2n-1}(R(\xi))
\end{equation}
where $a_{2n-1}(R(\xi))$ denotes the next to the smallest eigenvalue of $R(\xi)$.
\end{remark}

\begin{remark}
The choice of the nuclear norm and the square in (\ref{eq:l2bound}) is somewhat arbitrary. For any $p,q\geq 1$ (including infinity), there is
a constant $a_{p,q}>0$ so that 
\begin{equation}
\sum_{k=1}^m |\ip{F_k}{\outp{u}{v}}|^p \geq a_{p,q} \norm{\outp{u}{v}}_q^p.
\end{equation}
\end{remark}

An interesting corollary, which follows, is obtained in the case when $\nonlin$ is restricted to a subspace of $H$. 
It turns out that sometimes the underlying signal is actually real.
The canonical description of such a condition is to be conjugation invariant. 
Let $H'$ denote this set:
\begin{equation}
H'=\{x\in H~~;~~c(x)=x\}.
\end{equation}
Note $H'$ is a $\R$-linear space, but it is not a $\C$-linear space. Since $x$ is restricted to $H'$ it follows that the equivalence class in $H'$ of $x$ is
given by $\hat{x}\cap H'=\{ x,-x\}$. Consequently the appropriate quotient space is given by $\widehat{H'}=\{\{x,-x\}~,~x\in H_\R\}$.
Let $\pi_1:H\times H\rightarrow H$ and $\pi_2:H\times H\rightarrow H$ be the
canonical projections onto factors: $\pi_1((u,v))=u$ and $\pi_2((u,v))=v$. Then it is immediate to check that $H'$ admits the following equivalent descriptions:
\begin{equation}
 H' = \{ x\in H~;~\pi_2(\i(x))=0\} = \i^{-1}\left(\{H_\R\cap (H\times \{0\}) \right).
\end{equation}
Note in $H'$, $\ip{u}{v}$ is always real for all $u,v\in H'$ since $\ip{u}{v}=\ip{c(v)}{c(u)} = \ip{v}{u}$. 
Let $\fc=\{f_1,\ldots,f_m\}$ be the frame set in $H$. Note we do not assume $\fc\subset H'$.
Let
\begin{equation}
g_k = \pi_1(\i(f_k))~~,~~h_k=\pi_2(\i(f_k))
\end{equation}
with $1\leq k\leq m$. Note $f_k=\i^{-1}((g_k,h_k))=g_k+ih_k$ and $g_k,h_k\in H'$. Set:
\begin{equation}
\label{eq:Phik}
\Phit_k = g_kg_k^*+h_kh_k^* \in\SS^{2,0}(H')\subset \SS^{2,0}(H)~~,~~1\leq k\leq m.
\end{equation}
Note:
\begin{eqnarray}
\label{eq:Phi}
\ip{F_k}{\outp{x}{x}} & = & |\ip{f_k}{x}|^2 = \ip{\Phit_k x}{x} = \ip{\Phit_k}{\outp{x}{x}} \\
\label{eq:Phiuv}
\ip{F_k}{\outp{u}{v}} & = & real(\ip{u}{f_k}\ip{f_k}{v}) = \ip{\Phit_k u}{v} = \ip{\Phit_k}{\outp{u}{v}}
\end{eqnarray}
for all $x,u,v\in H'$. 
Thus the linear map $\Ac$ restricted to $Sym(H')$ can be thought of as taking inner products with a family of rank-2 nonnegative operators.
We have

\begin{corollary}\label{cor3.1}
Let $H$ be an $n$-dimensional complex Hilbert space with scalar product $\ip{}{}$ and conjugation $c:H\rightarrow H$. Let $H'=\{x\in H~;~c(x)=x\}$
be the maximal $c$-invariant set. The following are equivalent:
\begin{enumerate}
\item The restriction to $H'$ of the nonlinear map $\nonlin{\vert}_{H'}$ is injective on $\widehat{H'}$.
\item There is a constant $a_0>0$ so that $\forall u,v\in H'$ 
\begin{equation}
\label{eq:l2boundR}
\sum_{k=1}^m |\ip{\Phit_k u}{v}|^2 \geq a_0 \norm{u}^2 \norm{v}^2.
\end{equation}
\item For all $u\neq 0$,
\begin{equation}\label{eq:l2boundR2}
\dim_{\R} span_{\R}\{\ip{u}{g_k}g_k+\ip{u}{h_k}h_k~;~1\leq k\leq m\} = n.
\end{equation}
\item For any $u\neq 0$,
\begin{equation}
rank\left( \sum_{k=1}^m \Phit_k uu^* \Phit_k \right) = n.
\end{equation}
\item Let $R'(u)=\sum_{k=1}^m \Phit_k uu^* \Phit_k$. There is a constant $a_0>0$ so that $\forall u,v\in H'$,
\begin{equation}
\label{eq:Ibound2}
\ip{R'(u)v}{v} \geq a_0 \norm{u}^2 \norm{v}^2.
\end{equation}
\end{enumerate}
\end{corollary}
The proof is given in section \ref{subsec3.3}.

\subsection{Analyis of sets $\SS^{p,q}$ \label{subsec3.2}}

In addition to (\ref{eq:Spqdecomp}) the sets $\SS^{p,q}$ introduced in (\ref{eq:Spq})  have the following properties summarized in the following lemma:
\begin{lemma}
\label{l3.2.1}

\begin{enumerate}

\item For any $p_1\leq p_2$ and $q_1\leq q_2$, $\SS^{p_1,q_1}\subseteq \SS^{p_2,q_2}$

\item For any nonnegative integers $p,q$ the following disjoint decomposition holds true
\begin{equation}
\label{eq:Spqdecomp2}
\SS^{p,q} = \bigcup_{r=0}^p\bigcup_{s=0}^q \mathring{\SS}^{r,s}
\end{equation}
where by convention $\mathring{\SS}^{0,0}=\SS^{0,0}=\{0\}$, and $\mathring{\SS}^{p,q}=\emptyset$ for $p+q>n$.

\item For any nonnegative integers $p,q$,
\begin{equation}
\label{eq:Spqminu}
-\SS^{p,q}=\SS^{q,p}
\end{equation}

\item The mapping $(T,X)\mapsto TXT^*$ defines an action of $B(H)$ on $\SS^{p,q}$. 
Specifically for any $T\in B(H)$ and integers $p,q$,
\begin{equation}
T \SS^{p,q} T^* \subseteq \SS^{p,q}
\end{equation}
The inclusion becomes equality if $T$ is invertible.

\item $GL(H)$ acts transitively on $\mathring{\SS}^{p,q}$. 
Specifically for any $X,Y\in\mathring{\SS}^{p,q}$ there is an invertible $T\in GL(H)$ so that $Y=TXT^*$.

\item For any integers $p,q,r,s$,
\begin{equation}
\label{eq:Ssum}
\SS^{p,q} + \SS^{r,s} = \SS^{p,q} - \SS^{s,r} = \SS^{p+r,q+s}
\end{equation}

\end{enumerate}
\end{lemma}

{\bf Proof of Lemma \ref{l3.2.1}}

First three assertions are trivial. 
%Assertions (4) \& (5) state that signature and rank are the only invariants that quadratic forms over $\C$ have. 

(4) Fix an orthonormal basis in $H$, fix a $T\in B(H)$ and let $T=UDV$ be its singular value decomposition, 
where $U,V$ are unitary operators on $H$ and $D$ is a diagonal operator with non-negative entries. Let $X\in\SS^{p,q}$
and set $R(t)=U(tD+(1-t)I)V XV^*(tD+(1-t)I)U^*$ for every $0\leq t\leq 1$, where $I$ denotes the identity operator on $H$. Note $R(0)=UVXV^*U^* \in\SS^{p,q}$ and $R(1)=TXT^*$. For every $0\leq t <1$, the operator 
$U(tD+(1-t)I)V$ is invertible. Then by Sylvester's Law of Inertia (see Ex. 12.43 in Chapter 12 of \cite{Noble69}), 
for every $0\leq t<1$ the operator 
$R(t)$ has the same number of strictly positive eigenvalues and strictly negative eigenvalues as $X$ does. 
Since the spectrum is continuous with respect to matrix entries, it follows the number of strictly positive eigenvalues
cannot increase when passing to limit $t\rightarrow 1$. Same conclusion holds for the number of strictly 
negative eigenvalues. This shows $TXT^*=R(1)\in\SS^{p,q}$. 
Finally, when $T$ is invertible, $TXT^*$ has the same number of strictly positive (negative) eigenvalues as $X$ does. This shows $T\SS^{p,q}T^*=\SS^{p,q}$.

%Assertion (4) can be shown
%using the singular value decomposition of $T$ and the continuity of spectrum with respect to matrix entries. 

(5)  The conclusion follows again from Sylvester's Law of Intertia. Indeed fix an orthonormal basis in $H$ and let $T_1, T_2\in GL(H)$ be invertible transformations so that both $T_1 X T_1^*$ and $T_2 YT_2^*$ have the same matrix representations that are diagonal with $+1$ repeated $p$ times, $-1$ repeated $q$ times and $0$ repeated $n-p-q$ times. Thus $T_1 XT_1^*=T_2YT_2^*$ from where the conclusion follows with $T=T_2^{-1}T_1$. 

(6) One can see this statement as a special instance of the Witt decomposition theorem, a much more powerful tool in the theory of quadratic forms especially
in the case of vector spaces over fields of charactersitics other than 2, see \cite{Lam73}. However for the benefit of those who prefer a more direct proof, here is
a sketch of this result.
 Using spectral decomposition and rearranging the term, one can easily see that $\SS^{p+r,q+s}\subset \SS^{p,q}+\SS^{r,s}$.  
For the converse inclusion, we need to show that if $T,S\in Sym(H)$ are so that $T$ has at most p positive eigenvalues and q negative eigenvalues, and $S$
has at most $r$ positive eigenvalues and $s$ negative eigenvalues, then $T+S$ has at most $p+r$ positive eigenvalues, and $q+s$ negative eigenvalues.
Without loss of generality we can assume $T\in\mathring{\SS}^{p,q}$ and $S\in\mathring{\SS}^{r,s}$. Using spectral decompostions of $T$ and $S$
we obtain
\[ 
\begin{split}
T+S = & \left(\sum_{k=1}^p a_k A_k - \sum_{k=p+1}^{p+q} a_k A_k \right) + \left(\sum_{k=1}^r b_k B_k - \sum_{k=r+1}^{r+s} b_k B_k \right) = U-V ~,~ \\
 U & = \sum_{k=1}^p a_k A_k + \sum_{k=1}^r b_k B_k\geq 0~,~V=\sum_{k=p+1}^q a_k A_k+\sum_{k=r+1}^{r+s} b_k B_k\geq 0 ~,
\end{split}
\]
where $A_k,B_k$ are rank one orthogonal (spectral) projectors with $A_kA_j=0$ and $B_lB_h=0$ for all $k\neq j$ and $l\neq h$, and $a_k,b_k>0$. 
Thus $U\in\SS^{p+r,0}$ and $V\in\SS^{q+s,0}$. The claim now follows by induction provided we show that for any $R\in\mathring{\SS}^{a,b}$ and $E\in\SS^{1,0}$, a rank one
associated to vector $e\in H$ ($E=\ip{\cdot}{e}e$), then $E+R\in\SS^{a+1,b}$. 

Indeed this last assertion is shown as follows. Let $R=\sum_{k=1}^a c_k 
\ip{\cdot}{g_k}g_k - \sum_{k=a+1}^{a+b} c_k  \ip{\cdot}{g_k}g_k$ be its spectral decomposition. Two cases are treated distinctly.

Case 1: If $e\not\in R(H)$ then  $\Gamma'=\{g_1,\ldots,g_{a+b},e\}$ is linearly
independent. Let $\Gamma=\{\gamma_k~,~1\leq k\leq a+b+1\}$ be the (unique) biorthogonal system to $\Gamma'$. Let $T$ be an invertible operator that maps an orthonormal set $\Delta=\{\delta_1,\ldots,\delta_{a+b+1}\}$ into $\Gamma$, $T\delta_k = \gamma_k$, $1\leq k\leq a+b+1$, and maps the orthogonal space to $\Delta$ onto the orthogonal 
complement to $\Gamma$. Note $T^* g_k=\delta_k$, $1\leq k\leq a+b$ and $T^*e=\delta_{a+b+1}$. Then a direct computation shows that
\[ T^*(R+E)T = \sum_{k=1}^{a+b} c_k \ip{\cdot}{\delta_k}\delta_k + \ip{\cdot}{\delta_{a+b+1}}\delta_{a+b+1}. \]
Thus the spectrum of $T^*(R+E)T$ is composed of $\{c_1,\ldots,c_{a},-c_{a+1},\ldots,c_{a+b},1\}$ which shows that $T^*(R+E)T\in \SS^{a+1,b}$, and by
 (4), $R+E\in\SS^{a+1,b}$. 

Case 2: $e\in R(H)$. The rank of $R+E$ is less than or equal to the rank of $R$. Hence $R+E\in\SS^{a',b'}$ with $a'+b'=a+b$. Now by
the continuity of spectrum with respect to small perturbations, it follows that, for a small perturbation $e\mapsto e'\not\in R(H)$, $R+E'\in\SS^{a'',b''}$ with
$a''\geq a',b''\geq b'$. But the proof of case 1 shows $\SS^{a'',b''}\subset \SS^{a+1,b}$. Hence $b\geq b''\geq b'\geq b$, and $a+1\geq a''\geq a'\geq a$.
Thus $R+E\in\SS^{a+1,b}$. \hspace{1in} $\Box$
\vspace{5mm}

Next we analyze the space $\Soneone$ in more detail. The special factorization we obtain here can be extended to other classes $\SS^{p,q}$ but we do not 
plan to do so here. 
First set the following matrices:
\begin{equation}
\label{eq:H}
\HHH = \left[ \begin{array}{cc}
0 & 1 \\
1 & 0 \end{array} \right] ~~,~~D = \left[ \begin{array}{cc}
1 & 0 \\ 0 & -1 \end{array} \right]~~,~~
V = \frac{1}{\sqrt{2}} \left[ \begin{array}{cc}
1 & 1 \\ -1 & 1 \end{array} \right].
\end{equation}
In the next result we use a generalized unitary group $U(1,1;\HHH)$. Recall its definition.

{\bf Definition} The groups $U(1,1)$ and $U(1,1;\HHH)$ are defined by
\begin{equation}
\label{eq:U11}
U(1,1) = \{ A\in \C^{2\times 2}~,~A^*DA = D \}
\end{equation}
\begin{equation}
\label{eq:U11H}
U(1,1;\HHH) = \{ A\in \C^{2\times 2}~,~A^*\HHH A = \HHH \}
\end{equation}
These groups have been studied extensively in literature. See for instance \cite{Simon2004}, section 10.4. In particular the two groups above
are unitarily equivalent to each other, and the matrix $V$ provides such an equivalence:
\begin{equation}
\label{eq:Veq}
A\in U(1,1;\HHH) \Leftrightarrow B = VAV^* \in U(1,1).
\end{equation}
 The quadratic form $\omega(z_1,z_2) = |z_1|^2-|z_2|^2$ is invariant under the action of $U(1,1)$, whereas
$\phi(z_1,z_2)=\bar{z}_1z_2 + z_1\bar{z}_2$ is invariant under the action of $U(1,1;\HHH)$.

\begin{lemma} \label{l3.2.2a}
\mbox{}

\begin{enumerate}
\item $\SS^{1,1}=\SS^{1,0}-\SS^{1,0} = \SS^{1,0} + \SS^{0,1} .$

\item For any $T\in\SS^{1,1}$ there are $u,v\in H$ so that 
\begin{equation}
\label{eq:T11}
T=\frac{1}{2} (uv^*+vu^*) = \outp{u}{v}.
\end{equation}

If $T=a_1e_1e_1^*-a_2e_2e_2^*$ with $a_1,a_2\geq 0$
and $\ip{e_k}{e_j}=\delta_{k,j}$ is its spectral factorization then 
\begin{equation}
\label{eq:u0v0}
u_0 = ~~\sqrt{a_1}e_1+\sqrt{a_2}e_2,~~v_0 =\sqrt{a_1}e_1-\sqrt{a_2}e_2
\end{equation}
provides a particular factorization in (\ref{eq:T11}). 
%Conversely, for any $u,v\in H$, the operator $T=\outp{u}{v}$ is in $\Soneone$.
\end{enumerate}
\end{lemma}

\begin{lemma}\label{l3.2.2b}
\mbox{}

\begin{enumerate}

\item Let $T=\outp{u}{v}$. Then traces and spectrum $Sp(T)=\{a_{+},a_{-}\}$ are given by
\begin{eqnarray}
tr\{T\} & = & real(\ip{u}{v}) = \ip{u}{v}_{\R} \label{eq:traceT} \\
tr \{  T^2 \} & = & \frac{1}{4} ((\ip{u}{v})^2 + (\ip{v}{u})^2 + 2\norm{u}^2\norm{v}^2) \label{eq:traceT2} \\
 & = & \frac{1}{2}\left( \norm{u}^2\norm{v}^2 +
\ip{u}{v}_{\R}^2 - \ip{iu}{v}_{\R}^2 \right) \nonumber \\
a_{+} & = & \frac{1}{2}\left( \ip{u}{v}_{\R} +\sqrt{\norm{u}^2\norm{v}^2 - \ip{iu}{v}_{\R}^2 } \right) \geq 0 \label{a+} \\
a_{-} & = & \frac{1}{2}\left( \ip{u}{v}_{\R} -\sqrt{\norm{u}^2\norm{v}^2 - \ip{iu}{v}_{\R}^2 } \right) \leq 0 \label{a-}
\end{eqnarray}
The nuclear norm of $T$ is given by
\begin{equation}
\label{eq:nucnorm}
\norm{T}_1 = a_{+}+|a_{-}| = \sqrt{\norm{u}^2\norm{v}^2 - \ip{iu}{v}_{\R}^2}.
\end{equation}
Hence $T\in\SS^{1,1}$.
\item Let $T=\outp{u_0}{v_0}\in\mathring{\SS}^{1,1}$. Then any pair $(u,v)$ of vectors, with $u,v\in H$ so that $T=\outp{u}{v}$ is given by
\begin{equation}\label{eq:uv}
u=a_{11}u_0+a_{12}v_0~,~v=a_{21}u_0+a_{22}v_0
\end{equation}
for some matrix $A=(a_{k,l})_{1\leq k,l\leq 2}$ with $A\in U(1,1;\HHH)$. Conversely, for any matrix $A\in U(1,1;\HHH)$, $\outp{u}{v}=\outp{u_0}{v_0}$
where $(u,v)$ are given by (\ref{eq:uv}).

\end{enumerate}
\end{lemma}

\begin{lemma}\label{l3.2.2c}
\mbox{}

\begin{enumerate}
\item Let $T=xx^*-yy^*$ for some $x,y\in H$. Then $T\in\Soneone$ with spectrum $Sp(T)=\{b_{+},b_{-}\}$ and traces given by
\begin{eqnarray}
tr\{ T \} & = & \norm{x}^2-\norm{y}^2 \label{eq:trace:T'2} \\
tr\{ T^2 \} & = & \norm{x}^4 + \norm{y}^4 - 2 |\ip{x}{y}|^2 \label{eq:traceT'22} \\
b_{\pm} & = & \frac{1}{2} \left( \norm{x}^2 - \norm{y}^2 \right) \pm \frac{1}{2} \sqrt{(\norm{x}^2 + \norm{y}^2)^2 - 4|\ip{x}{y}|^2} \\
\norm{T}_1 & = &   \sqrt{(\norm{x}^2 + \norm{y}^2)^2 - 4|\ip{x}{y}|^2} \label{eq:1normT}
\end{eqnarray}

\item Let $T=xx^*-yy^*\in\Soneone$. Any pair of vectors $(x',y')$, with $x',y'\in H$ so that $T=x'(x')^*-y'(y')^*$ is related
to $(x,y)$ via
\begin{equation}
\label{eq:x'y'}
x' = b_{11}x+b_{12}y~,~y'=b_{21}x+b_{22}y
\end{equation}
for some matrix $B=(b_{ij})_{1\leq i,j\leq 2}$ with $B\in U(1,1)$. Conversely, for any matrix $B\in U(1,1)$, $x'(x')^*-y'(y')^*=xx^*-yy^*$,
 where $x',y'$ are given by (\ref{eq:x'y'}).
\end{enumerate}

\end{lemma}

\begin{remark}
1. The need for studying $\Soneone$ arose from the behavior of the IRLS algorithm described in \cite{Bal12a}.
However,  it quickly became apparent that the space $\Soneone$ and its factorization given by (\ref{eq:T11}) are crucial
for understanding the injectivity of the nonlinear map $\nonlin$, especially in light of Theorem \ref{th2.3} (4), equation (\ref{eq:kerAc}).

2. The choice in (\ref{eq:u0v0}) has the following two additional properties:
\begin{eqnarray}
\norm{u_0} & = &\norm{v_0}=\sqrt{a_1+a_2} = \sqrt{a_{+}-a_{-}} = \sqrt{\norm{T}_1 } \\
\ip{u_0}{v_0} & = & a_1-a_2 = a_{+} + a_{-} = tr\{T\} ~~{\rm (a~real~number!)}
\end{eqnarray}
where $\norm{T}_1$ represents the nuclear norm of $T$, and $a_1=a_{+}\geq 0$ and $a_2=-a_{-}\geq 0$ are its singular eigenvalues.
\end{remark}

{\bf Proof of Lemma \ref{l3.2.2a}}

(1) is a direct application of Lemma \ref{l3.2.1}(5).

(2) Since $\outp{}{}$ is $\R$-linear and $\outp{e_1}{e_2}=\outp{e_2}{e_1}$ we obtain
\[ \outp{u_0}{v_0} = a_1 e_1e_1^* + \sqrt{a_1a_2} \outp{e_1}{e_2} -\sqrt{a_1a_2} \outp{e_1}{e_2} - a_2 e_2e_2^*
 = a_1e_1e_1^*-a_2 e_2e_2^* = T .\hspace{0.5in}\Box\]
%For the converse, note that if $T=\outp{u}{v}$ then $rank(T) = dim\,span\{u,v\}\leq 2$. The claim follows once we establish (\ref{a+}) and (\ref{a-}) which we do below.
%\hspace{1in}$\Box$

{\bf Proof of Lemma \ref{l3.2.2b}}

(1) The equation (\ref{eq:traceT}) comes from the definition (\ref{eq:outp}) and the fact that $tr(vu^*) = \ip{v}{u}$. For $T^2$ compute first
\[ T^2 = \frac{1}{4} \left( \ip{v}{u}vu^* + \norm{u}^2 vv^* + \norm{v}^2 uu^* + \ip{u}{v} uv^* \right). \]
Then (\ref{eq:traceT2}) follows from this relation and $real((\ip{u}{v})^2) = (real(\ip{u}{v}))^2 - (imag(\ip{u}{v}))^2$. 
Finally, (\ref{a+}) and (\ref{a-}) come from solving:
\begin{equation} 
\label{eq:eigs}
\begin{array}{rcl}
 \mbox{$a_{+}+a_{-}$} & = & \mbox{$tr(T)$} \\
\mbox{$a_{+}^2 + a_{-}^2$} & = & \mbox{$tr(T^2)$}
\end{array}
\end{equation}
and observing
\[ a_{+}a_{-}=\frac{1}{4}\left((\ip{u}{v}_\R)^2+(\ip{iu}{v}_\R)^2 - \norm{u}^2\norm{v}^2 \right) = \frac{1}{4}(|\ip{u}{v}|^2- \norm{u}^2\norm{v}^2) \leq 0
.
\]

(2) A direct computation shows
\[ 
\begin{split}
\outp{u}{v}  & = \,\frac{1}{2} (\bar{a}_{11}a_{21}+\bar{a}_{21}a_{11})u_0u_0^* + \frac{1}{2}( \bar{a}_{11}a_{22}+\bar{a}_{21}a_{12} )v_0u_0^* \\
 & + \, \frac{1}{2}( \bar{a}_{12}a_{21}+\bar{a}_{22}a_{11} )u_0v_0^* +\frac{1}{2} (\bar{a}_{12}a_{22}+\bar{a}_{22}a_{12} )v_0v_0^* .
\end{split}
\]
Since $\{u_0,v_0\}$ are linearly independent, $\outp{u}{v}=\outp{u_0}{v_0}$ implies
\[  \bar{a}_{11}a_{21}+\bar{a}_{21}a_{11} = 0 ~~,~~\bar{a}_{11}a_{22}+\bar{a}_{21}a_{12} = 1 \]
\[ \bar{a}_{12}a_{21}+\bar{a}_{22}a_{11} = 1 ~~,~~\bar{a}_{12}a_{22}+\bar{a}_{22}a_{12} = 0 .\]
which corresponds to $A^*\HHH A=\HHH$. Hence $A\in U(1,1;\HHH)$. Conversely, if $A\in U(1,1;\HHH)$ then the above relations are satisfied which imply $\outp{u}{v}=\outp{u_0}{v_0}$.
\hspace{1in}$\Box$

{\bf Proof of Lemma \ref{l3.2.2c}}

Claims (1) and (2) are similar to claims in lemma \ref{l3.2.2b} and follow by direct computation.\hspace{1in} $\Box$

\vspace{5mm}

Topologically, $\Sonezero$ and $\Soneone$ are not differentiable manifolds. Instead they are algebraic varieties since they are given by
the zero loci of certain polynomials. We have the following result:
\begin{lemma}
\label{l3.2.3}

(1) The set $\mathring{\SS}^{1,0}$ is an analytic manifold in $B(H)$ of real dimension $2n-1$. As a real manifold, its tangent space at $X=x_0x_0^*$ is
given by 
\begin{equation}
\label{eq:TSonezero}
\TTT_X \mathring{\SS}^{1,0} = \{ \outp{x_0}{y}~,~y\in H\}.
\end{equation}
The $\R$-linear embedding $H\mapsto T_X \mathring{\SS}^{1,0}$ given by $y\mapsto \varphi_{x_0}(y)=\outp{x_0}{y}$ has null space given by
$\ker\,\varphi_x = \{iax_0~;~a\in\R\}$. 

(2) The set $\mathring{\SS}^{1,1}$ is an analytic manifold of real dimension $4n-4$. As a real manifold, its tangent space at $X=\outp{x_0}{y_0}$ is given by
\begin{equation}
\label{eq:TSoneone}
\TTT_X \mathring{\SS}^{1,1} = \{ \outp{x_0}{u} + \outp{y_0}{v}~,~u,v\in H \}.
\end{equation}
The $\R$-linear embedding $H\times H\mapsto \TTT_X \mathring{\SS}^{1,1}$ given by $(u,v)\mapsto\varphi_{x_0,y_0}(u,v) = \outp{x_0}{u}+\outp{y_0}{v}$
has null space given by $\ker\,\varphi_{x_0,y_0} = \{a(ix_0,0) + b(0,iy_0) + c(y_0,-x_0) + d(iy_0,ix_0)~,~a,b,c,d\in\R\}$. 
\end{lemma}
 
{\bf Proof of Lemma \ref{l3.2.3} }

Let $c_1,\ldots,c_n:Sym(H)\rightarrow\R$ be the coefficients of the characteristic polynomial:
\[ det (sI-T) = s^n + c_1(T)s^{n-1} +c_2(T) s^{n-2} +\cdots + c_n(T) . \]
with $c_1(T) = -tr(T)$ and $c_n(T)=(-1)^n det(T)$. Note that the $c_j$'s are polynomials.

(1) The manifold structure can be shown as follows. First note that 
$$S^{+}=\{ S\in Sym(H)~;~c_1(T) = -tr(S) < 0 \}$$ is an open subset of $Sym(H)$ and therefore a 
manifold of same real dimension as $Sym(H)$ (which is $n^2$). Next note
\[ \mathring{\SS}^{1,0} = c_2^{-1}(0)\cap\cdots\cap c_n^{-1}(0)\cap S^{+}. \]
Hence $\mathring{\SS}^{1,0}$ is an algebraic variety. 
Next we obtain that $\mathring{\SS}^{1,0}$ is a homogeneous space and hence a real analytic manifold.
Indeed by Lemma \ref{l3.2.1}(5), $GL(H)$ acts transitively on $\mathring{\SS}^{1,0}$. Therefore it is
sufficient to verify the stabilizer group is closed. Fix $\{e_1,e_2,\ldots,e_n\}$ an orthonormal basis in $H$ and consider the rank-1 operator $X=e_1e_1^*$. 
The stabilizer group for $X$ is given by invertible transformations $T$ so that 
$Te_1=ze_1$ with  $z\in\C$, $|z|=1$.  With respect to the fixed orthonormal basis, the stabilizer is represented by the group of matrices of the form:
\[ \HH_X^{1,0}=\left\{ \left[ \begin{array}{cc}
\mbox{$e^{i\theta}$} & \mbox{$v$} \\
0 & \mbox{$M$}
\end{array} \right] ~~,~~\theta\in[0,2\pi)~,~v\in\C^{1\times (n-1)}~,~M\in\C^{(n-1)\times (n-1)}\, ,\, det(M)\neq 0 \right\}. \]
One can easily verify that $\HH_X^{1,0}$ is a closed subset of $GL(n,\C)$, the Lie group of $n\times n$ invertible complex matrices.
%By Cartan's theorem, $H_X$ is a Lie subgroup of $GL(n,\C)$, the Lie group of $n\times n$ invertible complex matrices.
Thus $\mathring{\SS}^{1,0}$ is diffeomorphic to the analytic manifold $GL(n,\C)/\HH_X^{1,0}$.

Next we determine the tangent space.
Let $X=x_0x_0^*\in\mathring{\SS}^{1,0}$. We consider the set of all differentiable curves 
$$\Gcal= \left\{\gamma:I\rightarrow \mathring{\SS}^{1,0}~,~\gamma(0)=X, 0\in I\subset\R ~{\rm open~interval}\right\}, $$ 
passing through $X$. Then the tangent space to $\mathring{\SS}^{1,0}$ at $X$ is given by
\[ \TTT_X \mathring{\SS}^{1,0} = \left\{ \frac{d}{dt} \gamma(t) {\vert}_{t=0} ~;~\gamma\in\Gcal \right\}. \]
For each such curve, $\gamma:I\rightarrow \mathring{\SS}^{1,0}$ there is a unique differentiable curve $x:I\rightarrow H$ such that 
$\gamma(t)=x(t)(x(t))^*$ with $x(0)=x_0$.
In fact, locally, 
\[ x(t) = \frac{1}{\sqrt{\ip{\gamma(t)(x_0)}{x_0}}}\gamma(t)(x_0) \]
which shows $x(t)$ is differentiable. Then a direct computation shows
\[ \frac{d}{dt} \gamma(t){|}_{t=0} = \outp{\dot{x}(0)}{x(0)} + \outp{x(0)}{\dot{x}(0)} = 2 \outp{x_0}{\dot{x}(0)}, ~{\rm for~ any}~\dot{x}(0)\in H.  \]
Since $\dot{x}(0)$ can be chosen arbitrarily in $H$, it follows the tangent space to $\mathring{\SS}^{1,0}$ at $X$ is given by (\ref{eq:TSonezero}).
The real dimension of the $\R$-vector space $\TTT_X \mathring{\SS}^{1,0}$ is $2n-1$ once we notice the kernel of the $\R$-linear map $\varphi_{x_0}$ is one dimensional
and given by the real span of $ix_0$. 

(2) The algebraic variety structure is given by the intersection
\[ \mathring{\SS}^{1,1} = c_3^{-1}(0)\cap\cdots\cap c_n^{-1}(0)\cap S^{--} \]
where
\[ S^{--} = \{ S\in Sym(H) ~;~ c_2(T)<0 \}. \]
Next we obtain that $\mathring{\SS}^{1,1}$ is a homogeneous space and hence a real analytic manifold.
Indeed by Lemma \ref{l3.2.1}(5), $GL(H)$ acts transitively on $\mathring{\SS}^{1,1}$. Therefore it is
sufficient to verify the stabilizer group is closed. Fix $\{e_1,e_2,\ldots,e_n\}$ an orthonormal basis in $H$ and 
consider the rank-2 operator $X=e_1e_1^*-e_2e_2^*\in\mathring{\SS}^{1,1}$. 
The stabilizer group for $X$ is given by invertible transformations $T$ whose matrix representations are of the form
\[ \HH_X^{1,1}=\left\{ \left[ \begin{array}{cc}
\mbox{$R$} & \mbox{$v$} \\
0 & \mbox{$M$}
\end{array} \right] ~~,~~R\in U(1,1)~,~v\in\C^{2\times (n-2)}~,~M\in\C^{(n-2)\times (n-2)}\, ,\, det(M)\neq 0 \right\} \]
where $U(1,1)$ was introduced in (\ref{eq:U11}).
One can easily verify that $\HH_X^{1,1}$ is a closed subset of $GL(n,\C)$, the Lie group of $n\times n$ invertible complex matrices.
%By Cartan's theorem, $H_X$ is a Lie subgroup of $GL(n,\C)$, the Lie group of $n\times n$ invertible complex matrices.
Thus $\mathring{\SS}^{1,1}$ is diffeomorphic to the analytic manifold $GL(n,\C)/\HH_X^{1,1}$.

Next we determine the tangent space.
Fix $X\in \mathring{\SS}^{1,1}$, $X=\outp{x_0}{y_0}$ with $\{x_0,y_0\}$ linearly independent, and let 
$$\Gcal' = \{\gamma:I\rightarrow \mathring{\SS}^{1,1}~,~ 0\in I\subset\R~{\rm open~interval}~,~\gamma(0)=X\}$$
 be the set of differentiable curves passing through $X$. Then the tangent space to $\mathring{\SS}^{1,1}$ at $X$ is given by
\[ \TTT_X \mathring{\SS}^{1,1} = \left\{ \frac{d}{dt} \gamma(t) {\vert}_{t=0} ~;~\gamma\in\Gcal' \right\}. \]
By Lemma \ref{l3.2.1} we know $\gamma(t)=\outp{x(t)}{y(t)}$
for some $x:I\rightarrow H$ and $y:I\rightarrow H$. Note these functions are not unique. However we can choose them to be given by the spectral factorization
of $\gamma(t)$ via (\ref{eq:u0v0}). Furthermore a direct application of holomorphic functional calculus (see section 148, Decomposition Theorem in \cite{Nagy}) shows that $T\rightarrow P_1 \in \mathring{\SS}^{1,0}$
and $T\mapsto P_2\in \mathring{\SS}^{1,0}$ are analytic, where $T=P_1-P_2$ is its spectral decomposition with $P_1P_2=0$. Hence $t\in I\mapsto P_1(t)$ and $t\in I\mapsto P_2(t)$ 
are differentiable. The component functions $x(t)$ and $y(t)$
now uniquely defined by:
\begin{equation}
x(t) = \frac{1}{\sqrt{\ip{P_1(t)(x(0)) }{x(0)} }}P_1(t)(x(0))~~,~~y(t) = \frac{1}{\sqrt{\ip{P_2(t)(y(0))}{y(0)}}}P_2(t)(y(0))
\end{equation}
are differentiable. The derivative at $t=0$ is given by
\[ \frac{d}{dt}\gamma(t){|}_{t=0} = 2\outp{x(0)}{\dot{y}(0)} + 2 \outp{y(0)}{\dot{x}(0)} \]
And since $(x(0),y(0))$ and $(x_0,y_0)$ are related by a $U(1,1;H)$ matrix, and the fact that $(\dot{x}(0),\dot{y}(0))$ is arbitrary in $H\times H$, we obtain
the tangent space is given by (\ref{eq:TSoneone}), that is:
\[ \TTT_X \mathring{\SS}^{1,1} = \{ \varphi_{x_0,y_0}(u,v) = \outp{x_0}{u} + \outp{y_0}{v}~,~u,v\in H \}.\] 
A direct computation shows that $(ix_0,0)$, $(0,iy_0)$, $(y_0,-x_0)$ and $(iy_0,ix_0)$ are the only independent vectors in the 
null space of the $\R$-linear map $(u,v)\mapsto\varphi_{x_0,y_0}(u,v)$ which implies $\dim_\R \TTT_X \mathring{\SS}^{1,1}=4n-4$. Hence the real dimension of the manifold $\mathring{\SS}^{1,1}$ is $4n-4$.\hspace{1in}$\Box$

\begin{remark}
Compared to the complex projective manifold $\C \PP^{n-1}=\PP(\C^n)$, $\mathring{\SS}^{1,0}$ is diffeomorphic to $\C \PP^{n-1}\times \R^{+}$. The extra $\R^{+}$ component
comes from the fact that rank-1 operators in $\mathring{\SS}^{1,0}$ have arbitrary trace. This explains the real dimension of $\mathring{\SS}^{1,0}$:
$2\dim_{\R}\C \PP^n +1=2(n-1)+1=2n-1$. On the other hand, using spectral factorization, $\mathring{\SS}^{1,1}$ has real dimension given by 
$\dim_{\R}\mathring{\SS}^{1,0} + \dim_{\R}\mathring{\SS}^{1,0} -2 =2(2n-1)-2=4n-4$. The $-2$ term comes from the orthogonality of the two eigenvectors
of an operator in $\mathring{\SS}^{1,1}$. The $4n-4$ dimension estimate has been derived also heuristically in
 \cite{BCMN13} right after proof of Lemma 9 (\cite{Fritz}).
\end{remark}
\begin{remark}
As suggested by Bernhard Bodmann \cite{BB13}, it can be shown that the subset of projections inside $\mathring{\SS}^{1,0}$ is in fact a K\"{a}hler manifold
 (diffeomorphic to $\C \PP^{n-1}$). 
%A similar conclusion holds for $\mathring{\SS}^{1,1}$ since the set of projections inside $\mathring{\SS}^{1,1}$ is the %same as the set of projections inside $\mathring{\SS}^{1,0}$. 
%However neither of these two manifolds, $\mathring{\SS}^{1,0}$ and $\mathring{\SS}^{1,1}$, is K\"{a}hler.
However $\mathring{\SS}^{1,0}$ is not a K\"{a}hler manifold.
\end{remark}

\subsection{Analysis of the linear map $\tau$ \label{subsec3.21}}

We introduced earlier the $\R$-linear map $\tau$ that maps $Sym(H)$ operators into $Sym(H_\R)$ operators, using the real linear structure on these spaces.
In order for the diagram (\ref{eq:tauT}) to be commutative, the  map $\tau(T)$ is given explicitly by
\begin{equation}
\tau(T)(\xi) = \i(T(\i^{-1}(\xi)) ~~,~~\xi\in H_\R
\end{equation}
The following lemma summarizes the basic properties of the map $\tau$.
\begin{lemma}\label{l3.2.4}

\begin{enumerate}
\item Let $P$ be an orthogonal projection of rank $k$ in $H$. Then $\tau(P)$ is an orthogonal projection of rank $2k$ in $H_\R$. Furthermore if $\{e_1,\ldots,e_k\}$
is an orthonormal basis in the range of $P$, then $\{\i(e_1),\ldots,\i(e_k),J\i(e_1),\ldots,J\i(e_k)\}$ is an orthonormal basis in the range of $\tau(P)$.
\item If $T\in Sym(H)$ has spectrum $(a_1,a_2,\ldots,a_n)$ then $\tau(T)$ in $Sym(H_\R)$ has spectrum $(a_1,a_1,a_2,a_2,\ldots,a_n,a_n)$. 
\item For any two operators $T,S\in Sym(H)$, $\tau(T),\tau(S)\in Sym(H_\R)$ and
\begin{equation}
\label{eq:trH}
tr\{ \tau(T)\tau(S)\} = \ip{\tau(T)}{\tau(S)}_{B(H_\R)} = 2  \ip{T}{S}_{B(H)} = 2 tr\{TS\}
\end{equation}
\item Let $1\leq p\leq\infty$. The $p$-norms of a symmetric operator $T\in Sym(H)$ and $\tau(T)\in Sym(H_\R)$, are related by
\begin{eqnarray}
\label{eq:normp}
\norm{\tau(T)}_p & = & 2^{1/p} \norm{T}_p ~,~{\rm if~}p<\infty \\
\norm{T} & =&  \norm{\tau(T)}_{\infty} = \norm{T}_{\infty} = \norm{T}
\end{eqnarray}
\end{enumerate}
\end{lemma}

{\bf Proof of Lemma \ref{l3.2.4}}

(1) First we prove the statement for rank-1 projections. This comes from directly checking equation (\ref{eq:P}).
Thus $P=ee^*$ gets mapped into $\tau(P)=\epsi\epsi^*+J\epsi\epsi^*J^*$, where $\epsi=\i(e)$. 
Next if $\{e_1,\ldots,e_k\}$ is an orthonormal basis for the range of $P$ then
\[ P = \sum_{l=1}^k e_le_l^* \]
By $\R$-linearity, $\tau(P)$ has the form
\[ \tau(P) = \sum_{l=1}^k (\epsi_l\epsi_l^* + J\epsi_l\epsi_l^*J^*) \]
where $\epsi_l=\i(e_l)$, $1\leq l\leq k$.  Note
\[ \ip{J\epsi_l}{J\epsi_s} = \ip{\epsi_l}{\epsi_s} = \ip{e_l}{e_s}_\R = \delta_{l,s} ~,~ \ip{J\epsi_l}{\epsi_s} = \ip{ie_l}{e_s}_\R = 0 \]
Hence $\{\epsi_1,\ldots,\epsi_k,J\epsi_1,\ldots,J\epsi_k\}$ is an orthonormal set and since it is spanning the range of $\tau(P)$ 
it is an orthonormal basis in the range of $\tau(P)$. Hence $\tau(P)$ is an orthonormal projection on $H_\R$ of rank $2k$.

(2) Follows by using the spectral factorization of $T$,
\begin{equation}
\label{eq:eq:eq}
 T = \sum_{k=1}^r b_k P_k \mapsto \tau(T) = \sum_{k=1}^r b_k \tau(P_k) 
\end{equation}
where $P_1,\ldots,P_r$ are spectral projections and $b_1,\ldots,b_r$ are their associated distinct eigenvalues. 
For all $k\neq l$, $P_kP_l=0$ which implies $\tau(P_k)\tau(P_l) = 0$. 
Thus the right hand side of the second equation in (\ref{eq:eq:eq}) represents the spectral factorization of $\tau(T)$.
Hence each $b_k$ is an eigenvalue of $\tau(T)$ but with multiplicity twice the multiplicity as an eigenvalue of $T$.  
The conclusion now follows.

(3) It is enought to show $tr\{\tau(T)\tau(S)\}=2tr\{TS\}$. Fix an orthonormal basis in $H$, say $\{e_1,\ldots,e_n\}$.
Then by (1) $\{\i(e_1),\ldots,\i(e_n),J\i(e_1),\ldots,J\i(e_n)\}$ is an orthonormal basis in $H_\R$. Note $J\i(e_k)=\i(ie_k)$.
It follows
\begin{eqnarray*}
 tr\{\tau(T)\tau(S)\} & = & \sum_{k=1}^n \ip{\tau(S)\i(e_k)}{\tau(T)\i(e_k)} + \ip{\tau(S) J\i(e_k)}{\tau(T)J\i(e_k)} \\
 & = & \sum_{k=1}^n \ip{S e_k}{T e_k} + \ip{S ie_k}{T i e_k} = 2\sum_{k=1}^n \ip{Se_k}{T e_k} 
\end{eqnarray*}
which proves the claim.

(4)  Follows from (2): for finite $p$,
\[ \norm{\tau(T)}_p =\left( |a_1|^p +|a_1|^p + \cdots + |a_n|^p + |a_n|^p \right)^{1/p} = 2^{1/p} 
\left(|a_1|^p + \cdots+|a_n|^p \right)^{1/p} = 2^{1/p} \norm{T}_p \]
 whereas for $p=\infty$,
\[ \norm{\tau(T)}_{\infty} = \max \{|a_1|,|a_1|,\ldots,|a_n|,|a_n|\} = \max \{|a_1|,\ldots,|a_n|\} = \norm{T}_\infty. \]
$\Box$

\subsection{Proof of Theorem \ref{t3.1}\label{subsec3.3} and its Corollary \ref{cor3.1}}
\mbox{}
\vspace{3mm}

{\bf Proof of Theorem \ref{t3.1}}

(1)$\Leftrightarrow$(2). According to Theorem \ref{th2.3} (4), nonlinear map $\nonlin$ is injective iff $ker(\Ac)\cap(\Sonezero-\Sonezero) = \{0\}$.
But using Lemma \ref{l3.2.2a}(1) and (2) we get equivalently that $\nonlin$ is injective iff for all $u,v\in H$ with $\outp{u}{v}\neq 0$,
\[ \Ac(\outp{u}{v}) \neq 0 \]
Equivalently this means
\[ \sum_{k=1}^m |\ip{F_k}{\outp{u}{v}}|^2 >0 \]
Consider now the unit ball in $\Soneone$ with respect to the nuclear norm, say $S_1^{1,1}$. This set is compact in $Sym(H)$. Then let
\begin{equation}
a_0 = min_{T\in\Soneone,\norm{T}_1=1} \sum_{k=1}^m |\ip{F_k}{T}|^2 
\end{equation}
By homogeneity we obtain (\ref{eq:l2bound}). Then
\[ \ip{F_k}{\outp{u}{v}} = \frac{1}{2}(\ip{u}{f_k}\ip{f_k}{v}+\ip{v}{f_k}\ip{f_k}{u}) = real(\ip{u}{f_k}\ip{f_k}{v}) \]
and by Lemma \ref{l3.2.2b}(1),
\[ \norm{\outp{u}{v}}_1^2 = \norm{u}^2\norm{v}^2-\ip{iu}{v}_\R^2 = \norm{u}^2\norm{v}^2 - (imag(\ip{u}{v}))^2 \]
Putting together all previous derivations we obtain the equivalence $(1)\Leftrightarrow(2)$.

(2)$\Leftrightarrow$(4). Using Lemma \ref{l3.2.4} (3) we obtain (\ref{eq:l2bound}) is equivalent to
\begin{equation}
\label{eq:I2}
 \sum_{k=1}^m |\ip{\tau(F_k)}{\tau(\outp{u}{v})}|^2 \geq 4 a_0 \left[ \norm{u}^2\norm{v}^2 - (real(\ip{iu}{v}))^2 \right] 
\end{equation}
Now (\ref{eq:P}) and (\ref{eq:PQ}) imply
\[ \tau(F_k) = \outp{\varphi_k}{\varphi_k} + \outp{J\varphi_k}{J\varphi_k} = \outp{\varphi_k}{\varphi_k} + J\outp{\varphi_k}{\varphi_k} J^* \]
\[ \tau(\outp{u}{v}) = \outp{\xi}{\eta} + \outp{J\xi}{J\eta} = \outp{\xi}{\eta} +J\outp{\xi}{\eta}J^* \]
where $\varphi_k=\i(f_k)$ and $\xi=\i(u)$, $\eta=\i(v)$ and $J^*$ is the adjoint of $J$. A direct computation using $J^*=-J$ shows that
\[ \ip{\tau(F_k)}{\outp{J\xi}{J\eta}} = \ip{\tau(F_k)}{\outp{\xi}{\eta}} \]
Thus
\[ \ip{\tau(F_k)}{\tau(\outp{u}{v})} = 2 \ip{\varphi_k\varphi_k^* + J\varphi_k\varphi_k^*J^*}{\outp{\xi}{\eta}} = 
2[\ip{\xi}{\varphi_k}\ip{\varphi_k}{\eta} + \ip{\xi}{J\varphi_k}\ip{J\varphi_k}{\eta}] \]
With the equation (\ref{eq:Iexp}) we obtain
\[ \ip{\tau(F_k)}{\tau(\outp{u}{v})} = 2 \ip{\Phi_k\xi}{\eta} \Rightarrow \sum_{k=1}^m |\ip{\tau(F_k)}{\tau(\outp{u}{v})}|^2 = 4 \ip{R(\xi)\eta}{\eta} \]
Now the right-hand side of (\ref{eq:I2}) is processed as follows. Note $\norm{u}=\norm{\xi}$, $\norm{v}=\norm{\eta}$, and
\[ real(\ip{iu}{v}) = \ip{iu}{v}_\R = \ip{J\xi}{\eta} \]
Thus
\[ \norm{u}^2\norm{v}^2-(real(\ip{iu}{v}))^2 = \norm{\xi}^2\norm{\eta}^2 - (\ip{J\xi}{\eta})^2 = \ip{(\norm{\xi}^2 1 - J\xi\xi^*J^*)\eta}{\eta} \]
Substituting in (\ref{eq:I2}) we obtain (\ref{eq:Ibound}).

(3)$\Leftrightarrow$(4). Assume $rank(R(\xi))=2n-1$ for all $\xi\neq 0$. A direct computation shows that $R(\xi)(J\xi) = 0$. Hence $J\xi$ is the only 
independent vector in $ker(R(\xi))$. It follows there is an $a=a(\xi)>0$ so that
\[ R(\xi) \geq a(\xi) P_{J\xi}^{\perp} \]
Note the $a(\xi)$ represents the smallest nonzero eigenvalue of $R(\xi)$ which must be the $2n-1^{th}$.  
Since the eigenvalues of a matrix depend continuously with the matrix entries, it follows that $a(\xi)$ is a continuous
function on $\xi$. Let $a_0=\min_{\norm{\xi}=1}a(\xi)$. 
Since the minimum is achieved somewhere on the unit sphere, $a_0>0$.
Using the homogeneity of degree 2 of $R(\xi)$, we get $a(\xi)=\norm{\xi}^2a(\frac{\xi}{\norm{\xi}})\geq a_0 \norm{\xi}^2$
which proves (\ref{eq:Ibound}). Conversely, if (\ref{eq:Ibound}) holds true, then $R(\xi)$ has rank at least $2n-1$. Again since $J\xi$ is in the kernel
of $R(\xi)$, it follows that $R(\xi)$ must be of rank exactly $2n-1$.\hspace{1in}$\Box$ 
\vspace{5mm}

{\bf Proof of Corollary \ref{cor3.1}}

(1)$\Leftrightarrow$(2) follows from Theorem \ref{t3.1} (2) and equation (\ref{eq:Phiuv}) and the fact that $imag(\ip{u}{v})=0$ for all $u,v\in H'$.

(2)$\Leftrightarrow$(5) follows from the relation
\[ \sum_{k=1}^m |\ip{\Phit_k u}{v}|^2 = \ip{ \left( \sum_{k=1}^m \Phit_k uu^*\Phit_k \right)  v}{v} .\]

(4)$\Leftrightarrow$(5). Note first $\dim_\R H'=n$ since $H_\R=\i(H') \oplus \i(iH')$ is an orthogonal decomposition of the $2n$-dimensional real space $H_\R$ into two isomorphic subspaces. 
Hence $R'(u)$ is of rank-$n$ if and only if it is bounded below by a multiple of the identity restricted to $H'$.  
Thus (\ref{eq:Ibound2}) follows by the homogeneity of $R'(u)$ with respect to $\norm{u}$.

(2)$\Rightarrow$(3). For $u\neq 0$, (\ref{eq:l2boundR}) implies $\{\Phit_k u~,~1\leq k\leq m\}$ spans $H'$. This is equivalent with
 (\ref{eq:l2boundR2}).

(3)$\Rightarrow$(2). From (\ref{eq:l2boundR2}) we obtain that $\{\Phit_k u~,~1\leq k\leq m\}$ is a frame for $H'$.
Then (\ref{eq:l2boundR}) follows from the lower frame bound condition. \hspace{1in}$\Box$

\section{Performance Bounds on Reconstruction Algorithms\label{sec4}}

In this section we present two performance bounds applicable to any reconstruction algorithm.
One bound is deterministic and is based on the constants $a_0$ introduced in Theorem \ref{t3.1}.
The second bound represents the Cramer-Rao Lower Bound for the stochastic model (\ref{eq:model}).

\subsection{Lipschitz bounds of the inverse map\label{subsec4.1}}

Consider the nonlinear map $\nonlin:H\rightarrow\R^m$. We shall establish a deterministic performance bound for any inversion algorithm
in terms of the Lipschitz bounds of the map:
\begin{equation}
\label{eq:4.1}
 \Ac: \Sonezero \rightarrow \R^m~~,~~\Ac(xx^*)=\nonlin(x) 
\end{equation}
Specifically we want to bound from above and below the following expression:
\begin{equation}
U(x,y) = \frac{\norm{\Ac(xx^*) - \Ac(yy^*)}^2}{\norm{xx^*-yy^*}_1^2}
\end{equation}
Since $xx^*-yy^*=\outp{u}{v}\in\Soneone$ for some $u,v\in H$ it follows:
\begin{eqnarray}
\sup_{x,y\in H} U(x,y) & = & \sup_{u,v\in H} \frac{\sum_{k=1}^m |\ip{F_k}{\outp{u}{v}}|^2}{\norm{\outp{u}{v}}_1^2} \\
\inf_{x,y\in H} U(x,y) & = & \inf_{u,v\in H} \frac{\sum_{k=1}^m |\ip{F_k}{\outp{u}{v}}|^2}{\norm{\outp{u}{v}}_1^2} 
\end{eqnarray}
These ratios can be further processed as follows
\[ \frac{\sum_{k=1}^m |\ip{F_k}{\outp{u}{v}}|^2}{\norm{\outp{u}{v}}_1^2} = \frac{\ip{R(\xi)\eta}{\eta}}{\norm{\xi}^2 \ip{P_{J\xi}^{\perp}\eta}{\eta}} \]
where $\xi=\i(u)$ and $\eta=\i(v)$. Since $R(\xi)\eta = R(\xi)P_{J\xi}^{\perp}\eta$ if follows:
\[ \sup_{\xi,\eta\neq 0}  \frac{\ip{R(\xi)\eta}{\eta}}{\norm{\xi}^2 \ip{P_{J\xi}^{\perp}\eta}{\eta}} = \sup_{\xi\neq 0}\frac{\norm{R(\xi)}}{\norm{\xi}^2}
 = \max_{\xi\in H_\R,\norm{\xi}=1} \norm{R(\xi)} \]
and
\[ \inf_{\xi,\eta\neq 0}  \frac{\ip{R(\xi)\eta}{\eta}}{\norm{\xi}^2 \ip{P_{J\xi}^{\perp}\eta}{\eta}} = a_0^{opt}. \]
Note the constant $a_0^{opt}$ obtained above is the same as the one given in (\ref{eq:a0}).
Thus we proved:
\begin{thm}
\label{t4.1}
Assume the nonlinear map $\nonlin$ is injective. 
Then the map $\Ac:\Sonezero\rightarrow\R^m$ defined in (\ref{eq:4.1}) is bi-Lipschitz between $(\Sonezero,\norm{\cdot}_1)$ and $(\R^m,\norm{\cdot})$ with the Euclidian norm,
and it has the upper Lipschitz bound
\begin{equation}\label{B0}
B_0 = \sqrt{\max_{\xi\in H_\R, \norm{\xi}=1} \norm{R(\xi)} }
\end{equation}
and the lower Lipschitz bound 
\begin{equation}
\label{eq:A0}
A_0 = \sqrt{a_0^{opt}} = \sqrt{\min_{\xi\in H_\R,\norm{\xi}=1} a_{2n-1}(R(\xi))}
\end{equation}
Specifically
\begin{equation}
A_0 \norm{xx^*-yy^*}_1 \leq \norm{\Ac(xx^*) - \Ac(yy^*)} \leq B_0 \norm{xx^*-yy^*}_1
\end{equation}
for all $x,y\in H$.
\end{thm}

\subsection{The Cramer-Rao Lower Bound (CRLB)\label{subsec4.2}}

Consider now the noise model (\ref{eq:model}). We are interested in obtaining a lower bound for any unbiased estimator of $x$.
The derivation of the CRLB in this paper coincides with the one presented in \cite{BCMN13} (see Theorem 23 there).
 In turn this follows the recipe presented in \cite{Bal12a}.
We will just present the key steps of this derivation. Note that our derivation is canonical, that is basis independent.  

%Throughout this subsection and the next sections we assume $H=\C^n$ with the usual Hilbert space structure (scalar product and conjugation).
%Consequently $\outp{u}{v}=\frac{1}{2}(uv^*+vu^*)$ where $u^*$ denote the conjugate transpose of vector $u$. Then
%$\i(x)=(real(x),imag(x))$, $H_\R=\R^{2n}$  and $\outp{\eta}{\eta}=\eta\eta^T$ is a $2n\times 2n$ symmetric matrix where $\eta^T$ is the transpose of vector $\eta$.

Due to non-holomorphy of the nonlinear map $\nonlin$, the analysis is done in the realification space $H_\R$. We denote $\zeta=\i(x)$ for the signal $x\in H$.
The frame set is $\fc=\{f_1,\ldots,f_m\}$ and $F_k=f_kf_k^*\in \Sonezero(H)$ denotes the measurement operators. 
Recall our notation $\Phi_k=\tau(F_k)=\varphi_k\varphi_k^*+J\varphi_k\varphi_k^*J^*\in \SS^{2,0}(H_\R)$.
First we compute the Fisher information matrix associated to $\zeta$.
The likelihood for this problem is
\begin{equation}
p(y|\zeta) = \frac{1}{(2\pi)^{m/2} \sigma^m}exp\left(  -\frac{1}{2\sigma^2}\norm{y-\nonlin(x)}^2 \right)=\frac{1}{(2\pi)^{m/2}\sigma^m} exp\left(
-\frac{1}{2\sigma^2}\sum_{k=1}^m |y_k - \ip{\Phi_k \zeta}{\zeta}|^2 \right)
\end{equation}
where $\sigma^2$ is the noise variance.
The Fisher information matrix is given by (see \cite{Kay2010})
\begin{equation}
I(\zeta) = \E \left[ (\nabla_\zeta \log(p(y|\zeta)))(\nabla_\zeta \log(p(y|\zeta)))^T \right].
\end{equation}
The canonical form of this operator is
\begin{equation}
I(\zeta) = \E \left[ \outp{\nabla_\zeta \log(p(y|\zeta))}{\nabla_\zeta \log(p(y|\zeta))} \right].
\end{equation}
A little bit of algebra shows
\begin{equation}
\label{eq:Fisher}
I(\zeta) = \frac{4}{\sigma^2} \sum_{k=1}^m \Phi_k \zeta\zeta^* \Phi_k = \frac{4}{\sigma^2}R(\zeta).
\end{equation}
In general the covariance of any unbiased estimator is bounded below by the inverse of the Fisher information matrix (operator).
However in this case the Fisher information operator is not invertible. This fact simply expresses the statement that $x$ is not identifiable 
from the measurements $y=\nonlin(x)\in\R^m$ alone. As we know the nonlinear map $\nonlin$ is not injective on $H$ but instead it is injective on $\hat{H}$.
The nonuniqueness on $H$ is reflected in having a singular Fisher information matrix on $H_\R$. To solve this issue we need to fix the 
global phase factor. One solution is to fix a basis and decide that a particular component (say the last component) is real. Such an approach
was taken in \cite{BCMN13}. Here we propose a canonical solution to this normalization. An oracle provides us with a vector $z_0\in H$,  
so that $\ip{x}{z_0}>0$ is positive real. Assume $z_0$ is normalized $\norm{z_0}=1$.
Note there are two pieces of information that can be extracted from here: First the fact that $x$ is not
orthogonal  to $z_0$; in particular $x\neq 0$. Second, the global phase to recover $x$ from its rank-1 operator $xx^*$ is uniquely
determined by the fact that $imag(\ip{x}{z_0})=0$ and $real(\ip{x}{z_0})>0$.

Under this scenario we want to analyze the Fisher information operator obtained earlier. Let $\psi_0=\i(z_0)\in H_\R$. We know
\[ \ip{\zeta}{\psi_0} = real(\ip{x}{z_0}) >0 ~~,~~\ip{\zeta}{J\psi_0} = imag(\ip{x}{z_0}) =0 \]
with $\zeta=\i(x)$. Let $\Pi$ denote the orthogonal projection onto the complement of $J\psi_0$,
\begin{equation}
\Pi:H_\R\rightarrow E~~,~~\Pi = 1 - J\psi_0\psi_0^*J^* 
\end{equation}
where $E=\{J\psi_0\}^{\perp}$. 
%The oracle tells us $\psi_0$. 
Let $H_{z_0}$ denote the following closed set
\begin{equation}
\label{Hz0}
H_{z_0} = \{ \xi\in H_\R ~~,~~\ip{\xi}{\psi_0} \geq 0\,,\, \ip{\xi}{J\psi_0}=0  \} \subset E.
\end{equation}
Note $\zeta$ belongs to the relative interior of $H_{z_0}$.
The class of estimators for $\zeta$
should include only functions
\begin{equation}
\omega:\R^m \rightarrow H_{z_0}
\end{equation}
In this case the appropriate Fisher information operator should be
\begin{equation}\label{eq:Izeta}
\tilde{I}(\zeta) := \Pi I(\zeta) \Pi = \frac{4}{\sigma^2} \sum_{k=1}^m \Pi \Phi_k \zeta\zeta^* \Phi_k \Pi
\end{equation}

The following lemma proves that under the scenario described here, $\tilde{I}(\zeta)$ is invertible on $H_{z_0}$.
\begin{lemma}
\label{l4.2}
Assume $\nonlin$ is injective on $\hat{H}$ and $z_0\in H$ is so that $\ip{x}{z_0}>0$. Let $\zeta=\i(x)$.
 Then 
\begin{equation}\label{eq:l42eq}
\tilde{I}(\zeta) := \Pi I(\zeta) \Pi \geq \frac{4}{\sigma^2}a_0 |\ip{x}{z_0}|^2\, \Pi
\end{equation}
where $a_0=a_0^{opt}$ is the same lower bound introduced in Theorem \ref{t3.1} whose optimal value is given by (\ref{eq:a0}). Furthermore this bound is tight.
\end{lemma}
 {\bf Proof}

Using (\ref{eq:I}) the left-hand side of (\ref{eq:l42eq}) is $\tilde{I}(\zeta) = \frac{4}{\sigma^2} \Pi R(\zeta)\Pi$.
We know $R(\zeta) \geq a_0 \norm{\zeta}^2 P_{J\zeta}^{\perp}$ from Theorem \ref{t3.1} (4).
Therefore we only need to show
\[ \norm{\zeta}^2 \,\Pi\, P_{J\zeta}^{\perp}\,\Pi\, \geq\, |\ip{\zeta}{\psi_0}|^2 \,\Pi \]
where $\psi_0=\i(z_0)$, and the inequality is tight. 
Without loss of generality we can assume $\norm{\zeta}=1$ since all expressions are homogeneous in $\norm{\zeta}$. Then we need to show that for any $\xi\in E$, $\norm{\xi}=1$,
\[ \ip{P_{J\zeta}^{\perp} \xi}{\xi} \geq |\ip{\zeta}{\psi_0}|^2 \]
This follows from
\[ \inf_{\norm{\xi}=1,\xi\in E} 1- |\ip{\xi}{J\zeta}|^2  = 1- \max_{\norm{\xi}=1,\xi\in E} |\ip{\xi}{J\zeta}|^2 = 1-\left| \ip{\frac{\Pi J \zeta}{\norm{\Pi J\zeta}}}{J\zeta} \right|^2 = |\ip{\zeta}{\psi_0}|^2 \]
The last equality follows by direct computation from:
\[ \Pi J \zeta = J(\zeta-\ip{\zeta}{\psi_0}\psi_0)~~,~~\norm{\zeta-\ip{\zeta}{\psi_0}\psi_0}^2 = 1 - |\ip{\zeta}{\psi_0}|^2 .\]
Note also the inequality in (\ref{eq:l42eq}) is tight since the lower bound is achieved for $\zeta=\psi_0=argmin_{\xi\in H_\R,\norm{\xi}=1}a_{2n-1}(R(\xi))$, the optimizer in (\ref{eq:a0}).\hspace{1in}$\Box$

Thus we established that $\tilde{I}(\xi)$ is invertible on $H_{z_0}$. See also \cite{BCMN13}, Lemma 22, for a similar statement. 

Recall an estimator $\omega:\R^m\rightarrow H_{z_0}$ is said to be {\em unbiased} if 
\begin{equation} \E[\omega(y)|\zeta=\i(x)] = \zeta\end{equation}
We can now state the main result of this section:
\begin{thm}
\label{t4.2}
Assume the nonlinear map $\nonlin$ is injective and fix a vector $z_0\in H$. For any vector $x\in H$ with $\ip{x}{z_0}>0$  the covariance operator of any unbiased estimator 
$\omega: \R^m\rightarrow H_{z_0}$ of $x$ is bounded below by the Cramer-Rao Lower Bound (CRLB) given by
\begin{equation}
\label{eq:CRLB}
Cov[\omega(y)|\zeta=\i(x)] \geq \frac{\sigma^2}{4} \left(\sum_{k=1}^m \Pi\Phi_k \zeta\zeta^* \Phi_k\Pi \right)^{\dagger}
\end{equation}
where $\dagger$ denotes the pseudoinverse operator, $\zeta=\i(x)$ and $\Pi=1-J\psi_0\psi_0^*J^*$. In particular the Mean Square Error of $\omega$ is bounded below by
\begin{equation}
\label{eq:CRLB-MSE}
MSE(\omega) = \E[\norm{x-\omega(y)}^2|\zeta=\i(x)] \geq \frac{\sigma^2}{4} tr\left\{\left(\sum_{k=1}^m \Pi\Phi_k \zeta\zeta^* \Phi_k\Pi \right)^{\dagger} \right\}
\end{equation}
\end{thm}
{\bf Proof}

The key observation is that $H_{z_0}$ is a relatively open subset of the real linear space $E$, the orthogonal complement of $J\psi_0$ in $H_\R$, $E=\{J\psi_0\}^{\perp}\cap H_\R$.
Consider an orthonormal basis in $E$ of the form $\{e_1,\ldots,e_{2n-1}\}$. Thus $\{e_1,\ldots,e_{2n-1},J\psi_0\}$
is an orthonormal basis in $H_\R$. The (column vector) gradient with respect to $E$, $\nabla^E_\zeta$ has the form $\nabla^E_\zeta = \Pi\nabla_\zeta$ where $\nabla$ is the 
gradient with respect to the local coordinates in $H_\R$ and $\Pi$ is the orthogonal projection onto $E$. 
This shows the Fisher information matrix associated to the Additive White Gaussian Noise (AWGN) measurement 
process (\ref{eq:model}) with $\zeta\in H_{z_0}$ is $\tilde{I}(\zeta)$ given by (\ref{eq:Izeta}). 
Theorem 3.2 in \cite{Kay2010} implies the covariance matrix of $\omega$ is bounded below by the inverse of $\tilde{I}(\zeta)$ restricted to $E$. This implies (\ref{eq:CRLB}). Equation (\ref{eq:CRLB-MSE}) follows from
$MSE(\omega)=tr\{Cov[\omega(y)|\zeta=\i(x)]\}$ and (\ref{eq:CRLB}).\hspace{1in}$\Box$
\vspace{5mm}

Lemma \ref{l4.2} allows us to predict an upper bound for the MSE of any efficient estimator (that is an unbiased estimator that achieves the CRLB):
\begin{corollary}
\label{cor4.3}
Assume $\omega:\R^m\rightarrow H_{z_0}$ is an unbiased estimator that achieves the CRLB (\ref{eq:CRLB}). Then its Mean-Square Error is bounded above by
\begin{equation}
\label{eq:MSE-upper}
MSE(\omega) = \E[\norm{x-\omega(y)}^2|x] \leq \frac{(2n-1) \sigma^2}{4a_0^{opt}|\ip{x}{z_0}|^2}
\end{equation}
\end{corollary}
{\bf Proof}

If $\omega$ is unbiased and achieves the CRLB (in other words, if $\omega$ is efficient) then 
\[ MSE(\omega) = \frac{\sigma^2}{4} tr\left\{\left(\sum_{k=1}^m \Pi\Phi_k \zeta\zeta^* \Phi_k\Pi \right)^{\dagger} \right\} = 
tr \left\{ \left(\tilde{I}(\zeta) \right)^{\dagger} \right \} \]
Then (\ref{eq:MSE-upper}) follows from (\ref{eq:l42eq}) by noting that $tr\{\Pi\} = 2n-1$.\hspace{1in}$\Box$

\section{The Iterative Regularized Least-Square (IRLS) Algorithm\label{sec5}}

Consider the additive noise model in (\ref{eq:model}). Our data is the vector $y\in\R^m$. Our goal is to find
an $x\in H$ that minimizes $\norm{y-\nonlin(x)}$, where we use the Euclidian norm. Set
\begin{equation}
J_0(X) = \sum_{k=1}^m |y_k -\ip{Xf_k}{f_k}|^2~,~J_0:Sym(H)\rightarrow \R
\end{equation}
and notice $J_0 (xx^*) = \norm{y-\nonlin(x)}^2$. 
The least-square error minimizer represents the Maximum Likelihood Estimator (MLE) when the noise is Gaussian.
In this section we discuss an optimization algorithm for this criterion.

Consider now $J_0=\norm{y-\Ac(X)}^2$ where is $X$ is restricted to $\mathring{\SS}^{1,0}$, which is an analytic manifold. 
Consider a current point $X^{(t)}=x^{(t)}(x^{(t)})^*$ in an iterative process. 
Then a descent direction can be thought of as a vector in the tangent space to the manifold. According to Lemma \ref{l3.2.3} (1), the tangent space
at $X^{(t)}$ is given by operators of the form $\outp{x^{(t)}}{\delta}$. Since $X^{(t)}+\outp{x^{(t)}}{\delta}=\outp{x^{(t)}}{x^{(t)}+\delta}\in\Soneone$, one would need to project
$X^{(t)}+\outp{x^{(t)}}{\delta}$ back into $\Sonezero$ and choose direction $\delta$ that minimizes (or at least decreases) $J_0(P(X^{(t)}+\outp{x^{(t)}}{\delta}))$, where $P$ is the 
(nonlinear) projection in $Sym(H)$ onto $\Sonezero$.
However since $J_0$ is well defined on $\Soneone$ we choose to optimize $\delta$ without projecting back onto $\Sonezero$. Thus we obtain the iterative process:
\[ x^{(t+1)} = argmin_u J_0(\outp{x^{(t)}}{u}) \]
However this process is not robust to noise, the main reason being ill-conditioning and multiple local minima of $J_0$ on $\Sonezero$. Instead we choose
to regularize this process and thus to introduce a different optimization criterion.

Consider the following functional
\begin{eqnarray}\label{eq:JJ}
&& J:H\times H\times \R^{+}\times\R^{+}\rightarrow \R^{+} \\
J(u,v,\lambda,\mu) & = & \sum_{k=1}^m|y_k - \frac{1}{2}(\ip{u}{f_k}\ip{f_k}{v}+\ip{v}{f_k}\ip{f_k}{u})|^2+
\lambda\norm{u}^2 +\mu\norm{u-v}^2+\lambda\norm{v}^2. \nonumber
\end{eqnarray}
Our ultimate goal is to minimize $J_0(uu^*)=\norm{y-\nonlin(u)}^2=J(u,u,0,\mu)$ over $u$, for some (and hence any) value of $\mu\in\R^{+}$.
Our strategy is based on the following iterative process:
\begin{alg}
{\bf The Iterative Regularized Least-Square (IRLS) Algorithm}

Step 0. Initialize $x^{0}$ as the global optimal solution for a specific pair $(\lambda_0,\mu_0)$.

Step 1. Iterate:

\hspace{5mm}1.1 Solve for
\begin{equation}\label{eq:xt1}
x^{(t+1)} = argmin_u J(u,x^{(t)},\lambda_t,\mu_t)
\end{equation}

\hspace{5mm}1.2 Update $\lambda_{t+1},\mu_{t+1}$ according to a specific policy;

Step 2.  Stop when some tolerance level is achieved.
\label{alg1} 
\end{alg}

As we describe below the update (\ref{eq:xt1}) can be modified to achieve a more robust behavior (see (\ref{eq:xt+1})).

\subsection{Initialization\label{subsec5.1}}
Consider the regularized least-square problem:
\[ min_u J(u,u,\lambda,0) = min_u \norm{y-\nonlin(u)}^2 + 2\lambda\norm{u}^2 \]
Note the following relation
\begin{eqnarray}
J(u,u,\lambda,0) &=& \norm{y}^2 +2\lambda\norm{u}^2 - 2\sum_{k=1}^m y_k |\ip{u}{f_k}|^2 + \sum_{k=1}^m
|\ip{u}{f_k}|^4 \nonumber \\
& = & \norm{y}^2 + 2\ip{(\lambda I-Q)u}{u}+\sum_{k=1}^m |\ip{u}{f_k}|^4 \label{eq:J0}
\end{eqnarray}
where
\begin{equation}
\label{eq:Q}
 Q=\sum_{k=1}^m y_k f_kf_k^* = \sum_{k=1}^m y_k F_k
\end{equation}
For $\lambda>\norm{Q}$ the optimal solution is $u=0$. 
Note that if  $Q\leq 0$ (as a quadratic form) then the optimal solution of $min_u \norm{y-\nonlin(u)}^2$ is $u=0$.
Consequently in the following we assume the largest eigenvalue of $Q$ is positive. As $\lambda$ decreases
the optimizer remains small. Hence we can neglect the forth order term in $u$ in the expansion above and obtain:
\[ J(u,u,\lambda,0) \approx \norm{y}^2 + 2\ip{(\lambda I-Q)u}{u} \] 
Thus the critical value of $\lambda$ for which we may get a nonzero solution is $\lambda=max eig(Q)$, which is the
maximum eigenvalue of $Q$. Let us denote by $a_1$ this (positive) eigenvalue and $e_1$ its associated normalized eigenvector.
This suggests to initialize $\lambda=\rho a_1$ for some $0<\rho\leq 1$ and $x^{(0)}=\beta e_1$, for some nonzero scalar $\beta$.
Substituting into (\ref{eq:J0}) we obtain
\[ J(\beta e_1,\beta e_1,\rho a_1,0) = \norm{y}^2 - 2(1-\rho)a_1\beta^2 + (\sum_{k=1}^m |\ip{e_1}{f_k}|^4)\beta^4 \]
For fixed $\rho$, the minimum over $\beta$ is achieved at
\begin{equation}\label{eq:beta0}
\beta_0=\sqrt{\frac{(1-\rho)a_1}{\sum_{k=1}^m|\ip{e_1}{f_k}|^4}}~~,~~x^{(0)}=\beta e_1
\end{equation}
The parameter $\mu$ controls the step size at each iteration. The larger the value the smaller the step. On the other
hand, a small value of this parameter may produce an unstable behavior of the iterates. In our implementation we use
the same initial value for both $\lambda$ and $\mu$:
\begin{equation}
\label{eq:mu0}
\mu_0 = \lambda_0 = \rho a_1
\end{equation}
 
\subsection{Iterations\label{subsec5.2}}

Minimization (\ref{eq:xt1}) is performed in the space $H_\R$. Let $\zeta^{(t)} =\i(x^{(t)})$ and $\xi=\i(u)$. Then 
\begin{equation}
\label{eq:Jxit}
J(u,x^{(t)},\lambda_t,\mu_t) = \sum_{k=1}^m |\ip{\left(\Phi_k \zeta^{(t)}(\zeta^{(t)})^*\Phi_k\right) \xi}{\xi} - y_k|^2 + \lambda_t \norm{\xi}^2
+ \mu_t\norm{\xi-\zeta^{(t)}}^2 + \lambda_t \norm{\zeta^{(t)}}^2 
\end{equation}
Note the criterion is quadratic in $\xi$. The unique minimum is given by solving the linear equation:
\begin{equation}
\label{eq:eta}
\left(\sum_{k=1}^m \Phi_k \zeta^{(t)}(\zeta^{(t)})^* \Phi_k + (\lambda_t+\mu_t)1 \right) \zeta^{(t+1)} = \left( \sum_{k=1}^m y_k \Phi_k  + \mu_t 1 \right) \zeta^{(t)}
\end{equation}
for $\zeta^{(t+1)}$. 
In our implementations we decrease $(\lambda_t,\mu_t)$  but we limit $\mu_t$ to a minimum value. Thus our adaptation policy is
\begin{eqnarray}
\lambda_{t+1} & = & \gamma \lambda_t \label{eq:lambdat}\\
\mu_{t+1} & = & max(\gamma\mu_t,\mu^{min}) \label{eq:mut}
\end{eqnarray}
where $0<\gamma<1$ is the rate parameter.

\subsection{Stopping Criterion\label{subsec5.3}}

One approach is to repeat the iterations until $\lambda$ reaches a preset value $\lambda^{min}$. As proved later in this section, the error is 
linearly dependent on $\lambda$.

Alternatively, one can stop the iterations once the modeling error becomes comparable to the noise variance. Specifically, a stopping criterion
could be
\begin{equation}
\sum_{k=1}^m |y_k-|\ip{x^{(t)}}{f_k}|^2|^2 \leq \kappa m\sigma^2
\end{equation}
where $\kappa\geq 1$, for instance $\kappa=3$.  

\subsection{Convergence and Optimality}

Consider the following three functionals $J_1,J_2,J_3:Sym(H)\times\R^+\times \R^+\rightarrow\R$ defined by
\begin{eqnarray}
J_1(X,\lambda,\mu) & = & \sum_{k=1}^m |y_k - \ip{X}{F_k}|^2 + 2(\lambda+\mu) \norm{X}_1 - 2\mu\, tr\{X\} \\
J_2(X,\lambda,\mu) & = & \sum_{k=1}^m |y_k - \ip{X}{F_k}|^2  +2\lambda\, a_{max}(X) - (2\lambda+4\mu)\, a_{min}(X) \\
J_3(X,\lambda,\mu) & = & \sum_{k=1}^m |y_k - \ip{X}{F_k}|^2  +2\lambda\, \norm{X}_1 -4\mu\, a_{min}(X)
\end{eqnarray}
where $a_{max}(X)$  and $a_{min}(X)$ are the maximum and minimum eigenvalue of $X$, respectively.
We can prove the following result:

\begin{lemma}\label{t5.1}
\begin{enumerate}
\item When restricted to $\Soneone$ the three criteria coincide:
\begin{equation}
\label{eq:J123}
J_1(X,\lambda,\mu) = J_2(X,\lambda,\mu) = J_3(X,\lambda,\mu)~,~\forall X\in \Soneone,\lambda,\mu\in\R.
\end{equation}
\item On $Sym(H)$, the three criteria $J_1,J_2,J_3$ are convex.
\item The minimum value of $J_1,J_2,J_3$ on $\Soneone$ coincides with the minimum value of $J$ on $H\times H$:
\begin{equation}\label{minJ}
\min_{X\in\Soneone}J_1(X,\lambda,\mu) = \min_{X\in\Soneone}J_2(X,\lambda,\mu) = \min_{X\in\Soneone}J_3(X,\lambda,\mu) = \min_{u,v\in H} J(u,v,\lambda,\mu)
\end{equation}
for any $\lambda,\mu\geq 0$. Any minimizer $\hat{X}\in\Soneone$ for $J_1,J_2,J_3$ and $(\hat{u},\hat{v})$ for $J$ satisfy
\begin{equation}
\hat{X}=\outp{\hat{u}}{\hat{v}}~,~\norm{\hat{u}}=\norm{\hat{v}}~,~imag(\ip{\hat{u}}{\hat{v}})=0
\end{equation}
\item Restricted to $\Sonezero$ all four criteria coincide:
\begin{equation}
J(u,u,\lambda,\mu) = J_1(uu^*,\lambda,\mu) = J_2(uu^*,\lambda,\mu) = J_3(uu^*,\lambda,\mu) = \norm{y-\nonlin(u)}^2+2\lambda\norm{u}^2
\end{equation}
and are independent of $\mu$.
\end{enumerate}
\end{lemma}

{\bf Proof}

For (1), the quadratic error term is the same in all three criteria, whereas the regularization terms are equal to each other:
\begin{eqnarray*}
2(\lambda+\mu)\norm{X}_1 -2\mu tr\{X\} & = & 2(\lambda+\mu)(a_1+a_2)-2\mu(a_1-a_2) = 2\lambda a_1 +(2\lambda + 4\mu)a_2 \\
2\lambda a_{max}(X) - (2\lambda+4\mu) a_{min}(X) & = & 2\lambda a_1 + (2\lambda + 4\mu)a_2 \\
2\lambda \norm{X}_1 - 4\mu a_{min}(X) & = & 2\lambda(a_1+a_2) - 4\mu (-a_2) = 2\lambda a_1 +(2\lambda +4\mu)a_2
\end{eqnarray*}
where $X=a_1 e_1e_1^*-a_2 e_2e_2^*$ with $a_1,a_2\geq 0$ and $\{e_1,e_2\}$ orthonormal set.

For (2) notice that the following four functions defined on the real vector space $Sym(H)$ are convex: $X\mapsto |y_k - \ip{X}{F_k}|^2$ , $X\mapsto \norm{X}_1$, $X\mapsto -tr\{X\}$, $X\mapsto a_{max}(X)$, whereas $X\mapsto a_{min}(X)$ is concave. The last two statements 
are consequences of the Weyl's Inequality, Theorem III.2.1 in \cite{Bhat1997} with $i=j=1$ in (III.5), and $i-j=n$ in (III.6).

For (3) and (4) note first the following relation:
\begin{eqnarray}
J(u,v,\lambda,\mu)-J_1(\outp{u}{v},\lambda,\mu) & = & (\lambda+\mu)\left[(\norm{u}-\norm{v})^2 + 2\norm{u}\norm{v}
% - \nonumber \right.\\
% & &  \left. 
-2\sqrt{\norm{u}^2\norm{v}^2-(imag(\ip{u}{v}))^2}\right] \nonumber \\
 & \geq & (\lambda+\mu)(\norm{u}-\norm{v})^2 \geq 0 \label{eq:eq}
\end{eqnarray}
that follows from (\ref{eq:traceT}) and (\ref{eq:nucnorm}). Using part (1) we obtain
\[ \min_{X\in\Soneone}J_1(X,\lambda,\mu) = \min_{X\in\Soneone}J_2(X,\lambda,\mu) = \min_{X\in\Soneone}J_3(X,\lambda,\mu) 
\leq \min_{u,v\in H} J(u,v,\lambda,\mu).
 \]
Let $\hat{X}$ denote the optimizer and let 
$\hat{X}=a_1\hat{e}_1\hat{e}_1^*-a_2\hat{e}_2\hat{e}_2^*$ be its spectral decomposition with $a_1,a_2\geq 0$ and $\ip{\hat{e}_i}{\hat{e}_j}=\delta_{i,j}$, $1\leq i,j\leq 2$. 
Set $\hat{u}=\sqrt{a_1}\hat{e}_1+\sqrt{a_2}\hat{e}_2$ and $\hat{v}=\sqrt{a_1}\hat{e}_1-\sqrt{a_2}\hat{e}_2$.  Note $\hat{X}=\outp{\hat{u}}{\hat{v}}$ 
and $\norm{\hat{u}}=\norm{\hat{v}}$,
 $imag(\ip{\hat{u}}{\hat{v}})=0$. Then (\ref{eq:eq}) implies $J(\hat{u},\hat{v},\lambda,\mu)=J_1(\hat{X},\lambda,\mu)$ which
proves (\ref{minJ}). 
Furthermore, let $(\hat{u},\hat{v})$ be a minimizer for $\min_{u,v\in H} J(u,v,\lambda,\mu)$ which achieves also equality in (\ref{minJ}). 
By (\ref{eq:eq}) it follows that $\norm{\hat{u}}=\norm{\hat{v}}$ and $imag(\ip{\hat{u}}{\hat{v}})=0$ which proves statement (3).

(4) follows from the first equality in (\ref{eq:eq}) and (\ref{eq:J123}). \hspace{1in}$\Box$

\begin{remark}
The criterion $J_2$ shows that the two regularization terms $\lambda(\norm{u}^2+\norm{v}^2)$ and $\mu\norm{u-v}^2$ have different effects on the optimizer:
the larger the parameter $\mu$ the closer the lower eigenvalue is to zero; hence the closer the optimizer $\hat{X}$ is to a rank-1 operator; on the other hand, the 
larger the parameter $\lambda$ the larger the cost of the $\Sonezero$ component in $\hat{X}$; hence $\norm{\hat{X}}$ remains bounded. 
\end{remark}

\begin{remark}
In the optimization problem $\inf_{u,v\in H}J(u,v,\lambda,\mu)$ the constraint is convex but the criterion is not jointly convex in $(u,v)$. It is however
convex in each individual variable $u$ and $v$. On the other hand in the optimization problem $\inf_{X\in\Soneone} J_s(X,\lambda,\mu)$, $1\leq s\leq 3$,
the criterion is convex in $X$, but the underlying constraint $X\in\Soneone$ does not define a convex set. 
\end{remark}

The optimization procedure outlined in Algorithm \ref{alg1} describes the following mechanism. Consider a path $(x^{(t)})_{t\geq 0}$ in $H$.
Set
\[ X^{(t+1)} = \outp{x^{(t)}}{x^{(t+1)}} = \frac{1}{2}\left(x^{(t+1)}(x^{(t)})^* + x^{(t)}(x^{(t+1)})^* \right) \]
Thus a trajectory $(x^{(t)})_{t\geq 0}$ in $H$ is mappped into a trajectory $(X^{(t)})_{t\geq 0}$ in $\Soneone$. Then the algorithm chooses $X^{(t+1)}$ 
along a tangent direction to $\Soneone$ at $\outp{x^{(t)}}{x^{(t-1)}}$, namely a direction of the form $\outp{x^{(t)}}{u}$ for some $u\in H$ that is of maximum
descent for $J_s$, $1\leq s\leq 3$. Since the algorithm performance is completely characterized by the sequence $(X^{(t)})_{t\geq 0}$, and since for a given $X\in\Soneone$
the minimum of $J(u,v,\lambda,\mu)$ over $u,v\in H$ subject to $\outp{u}{v}=X$ is achieved at a pair $(u,v)$ so that $\norm{u}=\norm{v}$
and $imag\ip{u}{v}=0$, we chose to rescale the vector $x^{(t+1)}$ to achieve norm $\sqrt{\norm{\zeta^{(t)}}\norm{\zeta^{(t+1)}}}$. Thus:
\begin{equation}\label{eq:xt+1}
x^{(t+1)}=\sqrt{\frac{\norm{\zeta^{(t)}}}{\norm{\zeta^{(t+1)}}}}\i^{-1}(\zeta^{(t+1)})
\end{equation}

\subsection{Robustness to noise and the effect of regularization}

In this subsection we make a stronger assumption of injectivity:

{\em Assumption A:} For every $x\in H$ there is a unique $X\in\Soneone$ so that $\Ac(X)=\Ac(xx^*)$, namely $X=xx^*$. 

A simple heuristic argument (that can be made more precise, but we do not intend to do so here) suggests that generically this assumption is satisfied
for frames of redundancy 6 or more (that is for $m\geq 6n$). 

This assumption turns out to be equivalent to a stability bound that we describe next

\begin{lemma} The following are equivalent:
\begin{enumerate}
\item The frame $\fc$ satisfies assumption A.
\item $\ker\Ac\cap\SS^{2,1} = \{0\}$.
\item There is a constant $A_3>0$ so that
\begin{equation}
\label{eq:XY}
A_3 \norm{X-Y}_1^2 \leq \norm{\Ac(X)-\Ac(Y)}^2
\end{equation}
for all $X\in\Soneone$ and $Y\in\Sonezero$.
\end{enumerate}
\end{lemma}

{\bf Proof}

(1)$\Rightarrow$(2). Note first that $\Sonezero-\Soneone =\SS^{2,1}$ cf. Lemma \ref{l3.2.1} (5). Assume $\fc$ satisfies assumption A
and let $Z\in\ker\Ac\cap\SS^{2,1}$. Then $Z=X-xx^*$ for some $X\in\Soneone$ and $x\in H$. Then $\Ac(X)=\Ac(xx^*)$ and by
Assumption A, $X=xx^*$. Hence $Z=0$.

(2)$\Rightarrow$(1).
Conversely if (2) holds true, then for every $x\in H$ and $X\in\Soneone$ so that $\Ac(X)=\Ac(xx^*)$ it follows $Z=X-xx^*\in\ker\Ac$
 and $Z\in\SS^{2,1}$. Thus $Z=0$ which means $\fc$ satisfies assumption A.

(3)$\Rightarrow$(2) is immediate.

(2)$\Rightarrow$(3). Since $S=\{ W\in\SS^{2,1}~|~\norm{W}_1=1\}$ is compact, it follows
\begin{equation}
\label{eq:A3}
A_3 = \inf_{W\in S} \norm{\Ac(W)}^2 = \norm{\Ac(W_0)}^2 >0
\end{equation}
for some $W_0\in\SS^{2,1}$, $W_0\neq 0$. Then since any $Z\in\SS^{2,1}$ can be written as $Z=\norm{Z}_1 W$, with $W=\frac{1}{\norm{Z}_1}Z\in S$,
\[ \norm{\Ac(Z)}^2 = \norm{Z}_1^2 \norm{\Ac(W)}^2 \geq A_3 \norm{Z}_1^2 \]
Then (\ref{eq:XY}) follows from noticing that $Y-X\in\SS^{2,1}$.\hspace{1in}$\Box$

Now we can state the main stability result of the estimators described in this section.

\begin{thm}
\label{t5.10}
Fix $\mu\geq 0$, $x\in H$ and let $y=\nonlin(x)+\nu$. Assume $\fc$ satisfies assumption A and let $A_3$ denote the Lipschitz bound in (\ref{eq:XY}).
Assume the optimization procedure finds a pair $u,v\in H$ so that $J(u,v,\lambda,\mu)\leq J(x,x,\lambda,\mu)$. Then
\begin{equation}
\label{eq:robust}
\norm{\outp{u}{v}-xx^*}_1 \leq \frac{2\lambda}{A_3}+2\sqrt{\frac{\lambda^2}{A_3^2} + \frac{\norm{\nu}^2}{{A_3}}} 
\leq 4\frac{\lambda}{A_3} + 2\frac{\norm{\nu}}{\sqrt{A_3}}~~.
\end{equation}
Let $\outp{u}{v}=a_1e_1e_1^*-a_2e_2e_2^*$, with real $a_1,a_2\geq 0$ and $\{e_1,e_2\}$ orthonormal set in $H$, be its spectral decomposition.
Assume an oracle provides the global phase $\varphi_0$ so that $e^{i\varphi_0}=\frac{\ip{x}{e_1}}{|\ip{x}{e_1}|}$. Set  
\begin{equation}
\hat{x}= e^{i\varphi_0} \sqrt{a_1}\,e_1~~.
\end{equation}
Then
\begin{equation}
\label{eq:l2err1}
\norm{x-\hat{x}}^2 \leq \norm{\outp{u}{v}-xx^*}_1 + a_2\leq \frac{4\lambda}{A_3}+\frac{2\norm{\nu}}{\sqrt{A_3}} + \frac{\norm{\nu}^2}{4\mu} + \frac{\lambda\norm{x}^2}{2\mu}~~.
\end{equation}
If, additionally, $J(u,v,\lambda,\mu)\leq J(0,0,\lambda,\mu)$ then
\begin{equation}\label{eq:l2err}
\norm{x-\hat{x}}^2 \leq \frac{4\lambda}{A_3} + \frac{2\norm{\nu}}{\sqrt{A_3}} + \frac{\norm{\nu}^2}{4\mu}~~.
\end{equation}
\end{thm}

{\bf Proof}

The proof follows from the Lipschitz bound (\ref{eq:XY}). Let $Y=\outp{u}{v}$ and $X=xx^*$. Note that
\[ J_3(Y,\lambda,\mu) \leq J(u,v,\lambda,\mu) \leq J(x,x,\lambda,\mu) = J_3(xx^*,\lambda,\mu) \]
where the first inequality follows from (\ref{eq:J123}) and (\ref{eq:eq}).
Explicitly this means
\[ \norm{\Ac(Y-X)-\nu}^2 + 2\lambda \norm{Y}_1 - 4\mu a_{min}(Y) \leq \norm{\nu}^2 + 2\lambda \norm{X}_1 ~~.\]
Then
\[ A_3 \norm{Y-X}_1^2 \leq \norm{\Ac(Y-X)}^2 \leq \left( \norm{\Ac(Y-X) -\nu} + \norm{\nu} \right)^2 \leq
2 \norm{\Ac(Y-X)-\nu}^2 + 2\norm{\nu}^2 ~~.\]
Substituting into the previous inequality we obtain
\[ \norm{X-Y}_1^2 -\frac{4\lambda}{A_3} \norm{X-Y}_1 - \frac{4\norm{\nu}^2  +8\mu a_{min}(Y)}{A_3} \leq 0 ~~.\]
Solving for $\norm{X-Y}_1$ we obtain 
\[ \norm{X-Y}_1 \leq \frac{2\lambda}{A_3} + 2\sqrt{\frac{\lambda^2}{A_3^2}+\frac{\norm{\nu}^2}{A_3}  + \frac{2 \mu}{A_3}a_{min}(Y)}
\leq \frac{2\lambda}{A_3} + 2\sqrt{\frac{\lambda^2}{A_3^2}+\frac{\norm{\nu}^2}{A_3}} \]
since $\mu\geq 0$ and $a_{min}(Y)\leq 0$ for any $Y\in\Soneone$. This proves the first inequality in (\ref{eq:robust}). 
The second inequality in (\ref{eq:robust}) follows from $\sqrt{a^2+b^2}\leq a+b$ for any $a,b\geq 0$. 

The second part of the Theorem is obtained as follows. 
First note $\hat{x}\hat{x}^* = a_1 e_1e_1^*$. Hence
\[ \norm{\hat{x}\hat{x}^* - xx^*}_1 \leq \norm{\outp{u}{v}-xx^*}_1 + a_2 ~~.\]
Next we show that $\norm{x-\hat{x}}^2\leq\norm{\hat{x}\hat{x}^*-xx^*}_1$, from where
the first inequality of (\ref{eq:l2err1}) follows.
Let $T=\hat{x}\hat{x}^*-xx^*\in\Soneone$. Its nuclear norm is given by (\ref{eq:1normT}):
\[ \norm{T}_1 =\sqrt{(\norm{\hat{x}}^2+\norm{x}^2)^2 - 4|\ip{x}{\hat{x}}|^2}~~. \]
Recall $\hat{x}$ is given the global phase so that $\ip{\hat{x}}{x}\geq 0$ is real and nonnegative. Thus $|\ip{\hat{x}}{x}|^2 = (\ip{\hat{x}}{x})^2$ and
\[ \norm{x-\hat{x}}^4 = \norm{T}_1^2 -4 \ip{\hat{x}}{x}\cdot\norm{x-\hat{x}}^2 \leq \norm{T}_1^2 \]
which shows $\norm{x-\hat{x}}^2\leq\norm{\hat{x}\hat{x}^*-xx^*}_1$ and thus first inequality in (\ref{eq:l2err1}). 
The second inequality follows from  (\ref{eq:robust}) and a bound for $a_2$ from
$$4\mu a_{2}\leq J_3(Y,\lambda,\mu)\leq J_3(X,\lambda,\mu)=\norm{\nu}^2+2\lambda\norm{x}^2.$$

Finally, (\ref{eq:l2err}) follows from using the second inequality in (\ref{eq:robust}) together with a bound for $a_{2}$ from 
$$4\mu|a_2|\leq J_3(Y,\lambda,\mu)\leq J_3(0,\lambda,\mu)=\norm{\nu}^2.$$
$\Box$

%\begin{remark}
%The assumption $J(u,v,\lambda,\mu)\leq J(0,0,\lambda,\mu)$ is not critical. A different upper bound than %(\ref{eq:l2err}) can  be obtained from 
%$-4\mu a_{2}(X)\leq J_3(Y,\lambda,\mu)=\norm{\nu}^2+2\lambda\norm{x}^2$:
%\begin{equation}
%\label{eq:err2}
%\norm{x-\hat{x}}^2 \leq 6\frac{\lambda}{A_3} + 4\norm{\nu} + \frac{\norm{\nu}^2}{4\mu} + %\frac{\lambda\norm{x}^2}{2\mu}
%\end{equation}
%which can be made small (to the order of noise level) by choosing a large $\mu$ and small $\lambda$. 
%\end{remark}

\begin{remark}
Note the result does not require to finding the global minimum, but only an estimate that brings the criterion $J$ below the level achieved by the true signal. 
On the other hand the result is not unexpected since MLE is asymptotically efficient (i.e. unbiased and achieves CRLB asymptotically
with the number of measurements), see Theorem 7.1 in \cite{Kay2010}.
\end{remark}

\section{Numerical Analysis\label{sec6}}

In this section we present numerical simulations for the Regularized Iterative Least-Square algorithm presented in the previous section.

We generated random frames or redundancy $\frac{m}{n}=$ 4,6, and 8, as well as complex random signals $x$
of size $n=100$.
All these vectors (frame and signal) are drawn from standard Gaussian distribution $\mathcal{N}(0,I)$. Then we scale the frame vectors to have norm 1.
We set the first component of $x$ to be a  positive real, and so the global phase becomes uniquely determined. 
The algorithm stopped when $\lambda_t$ reached a preset value. In these simulations we choose
$\lambda^{min}=0.01$, $\mu^{min}=1$ and rate $\gamma=0.8$ in (\ref{eq:lambdat},\ref{eq:mut}).

The magnitude square of signal coefficients $\nonlin(x)$ is perturbed additively by Gaussian noise with variance $\sigma^2$ 
to achieve a fixed signal-to-noise-ratio defined as
\[ SNR = \frac{\sum_{k=1}^m|\ip{x}{f_k}|^4}{m\sigma^2} ~~,~~SNRdB = 10\log_{10} (SNR)~~{\rm [dB]}\]

We vary SNRdB over 15 values in 5dB increments from -30dB to +40dB. We average algorithm performance over 100 noise realizations.

\begin{figure}[htb]
\includegraphics[width=130mm,height=70mm]{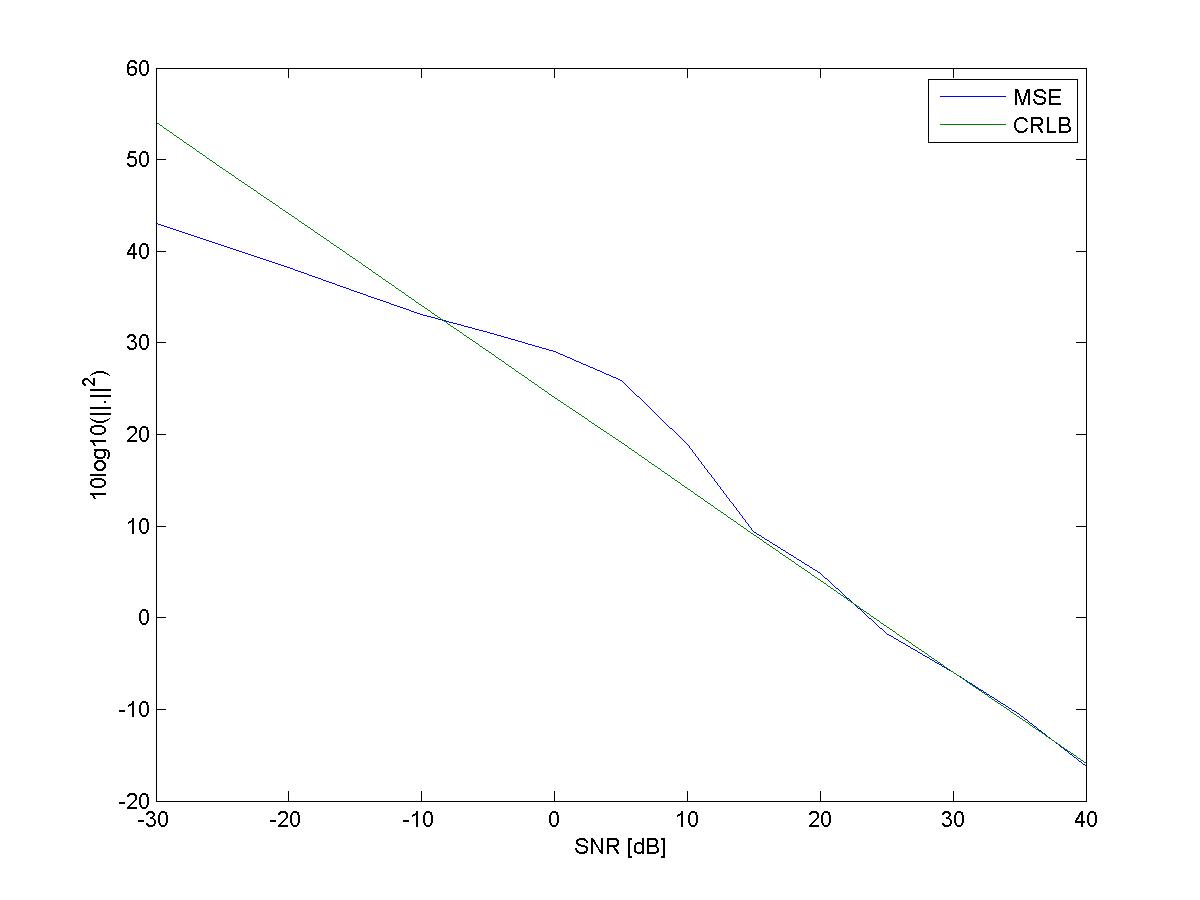}
\includegraphics[width=130mm,height=70mm]{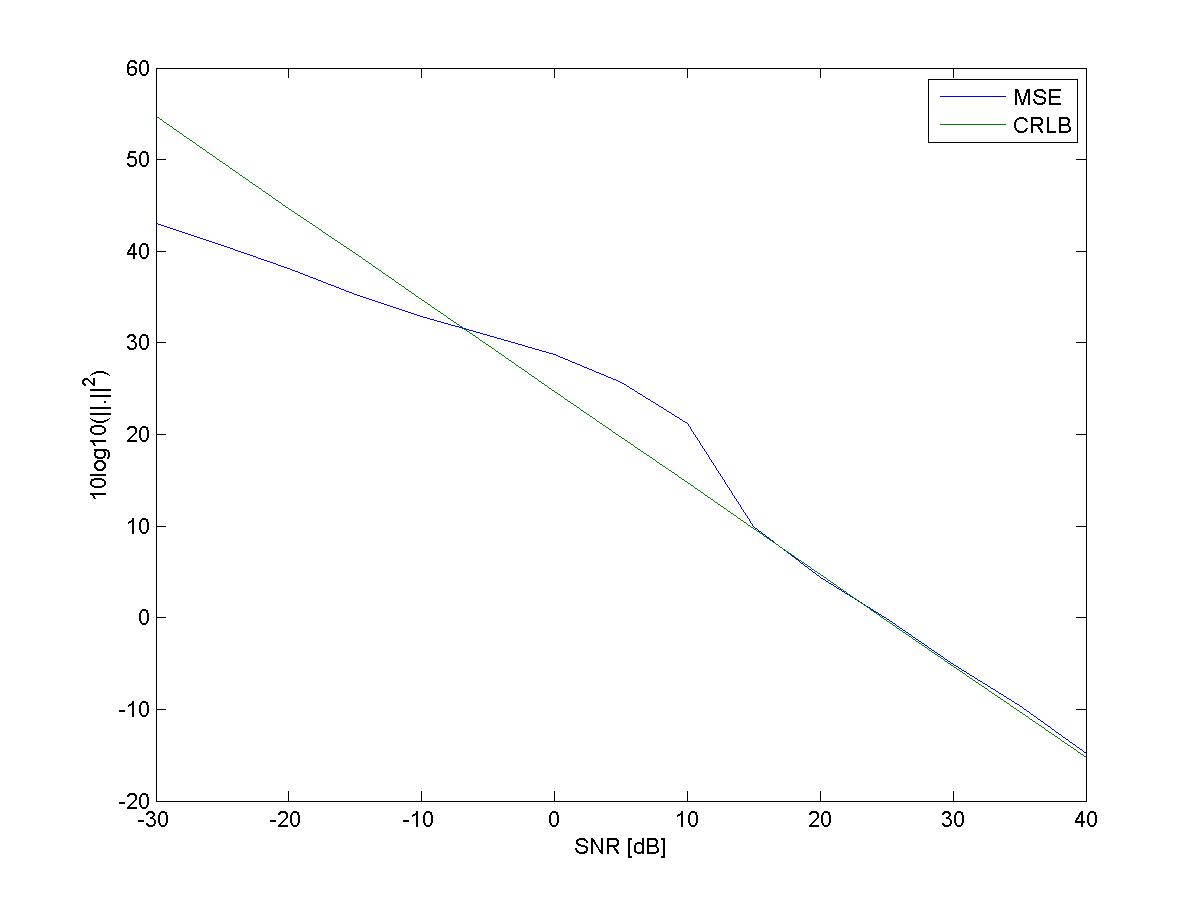}
\includegraphics[width=130mm,height=70mm]{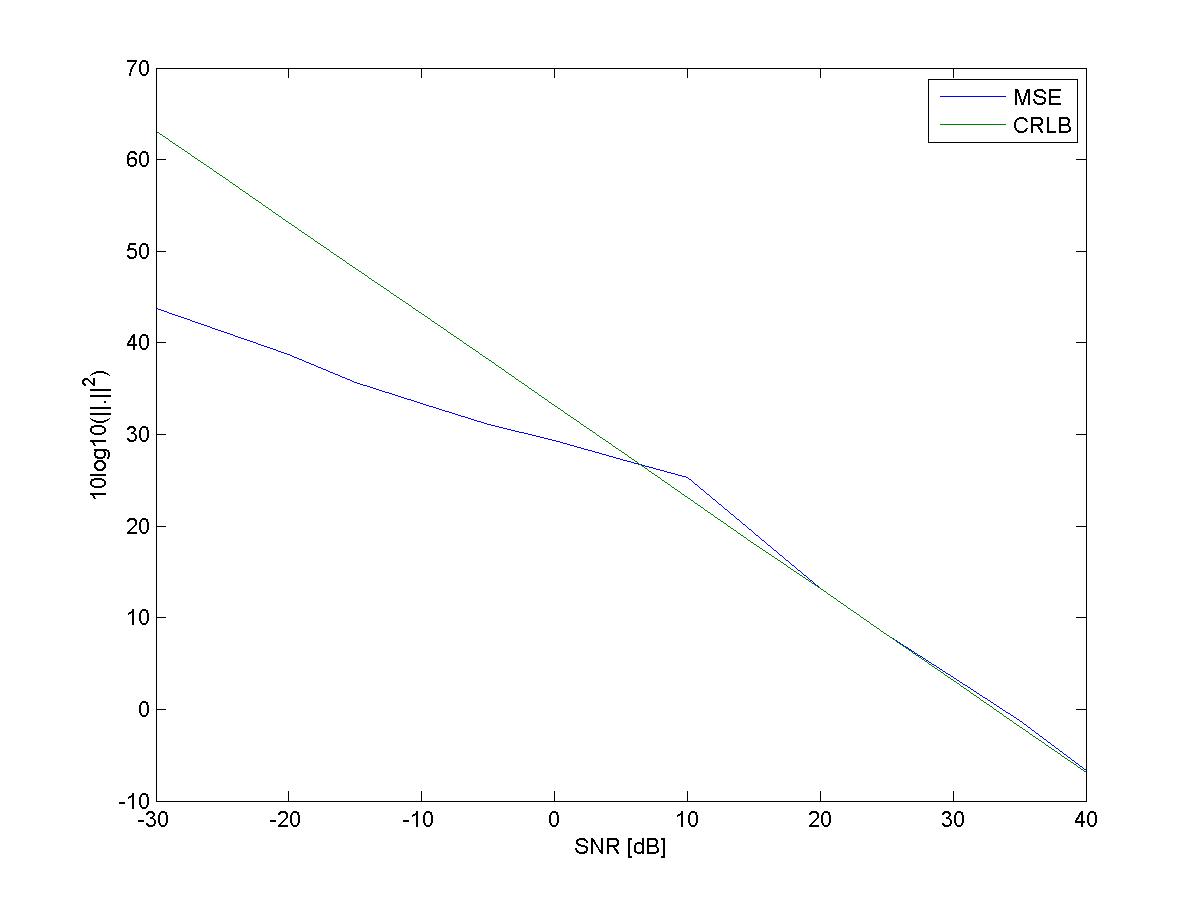}
\caption{Mean-Square Error and CRLB bounds for $m=800$ (top plot), $m=600$ (middle plot),
 and $m=400$ (bottom plot) when $n=100$. \label{fig1}}
\end{figure}

Figure \ref{fig1} includes the mean-square error averaged over 100 noise realizations for one fixed realization of $x$, and the C-R lower
bound. 

\begin{figure}[htb]
\includegraphics[width=130mm,height=70mm]{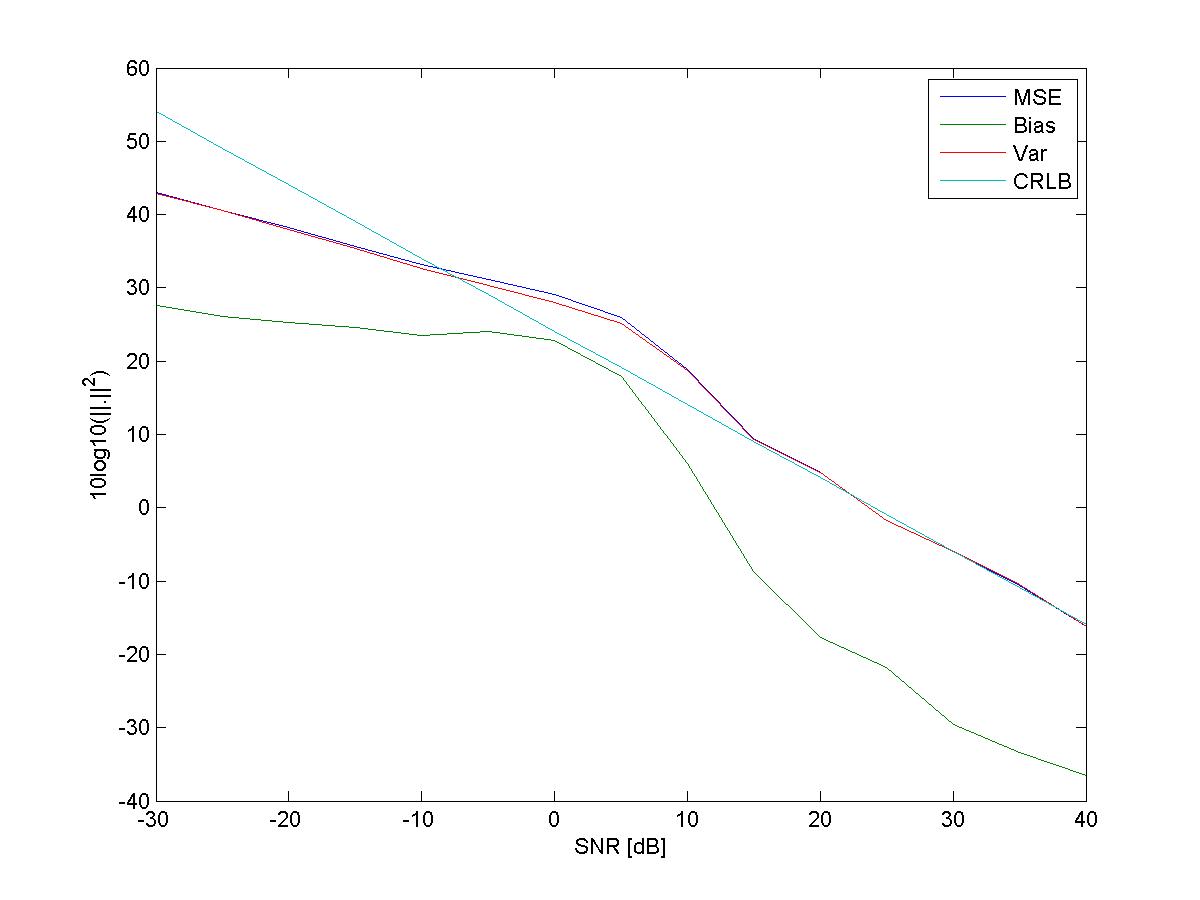}
\includegraphics[width=130mm,height=70mm]{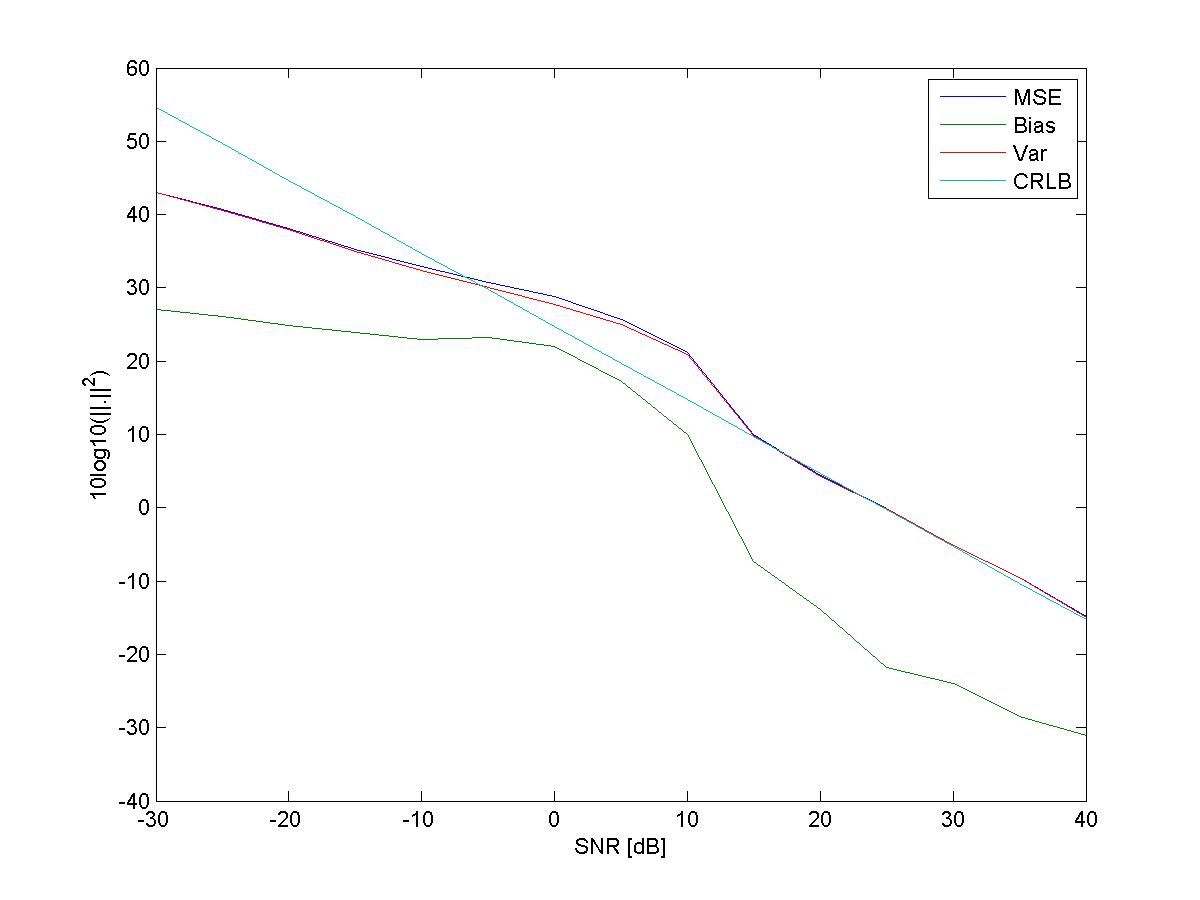}
\includegraphics[width=130mm,height=70mm]{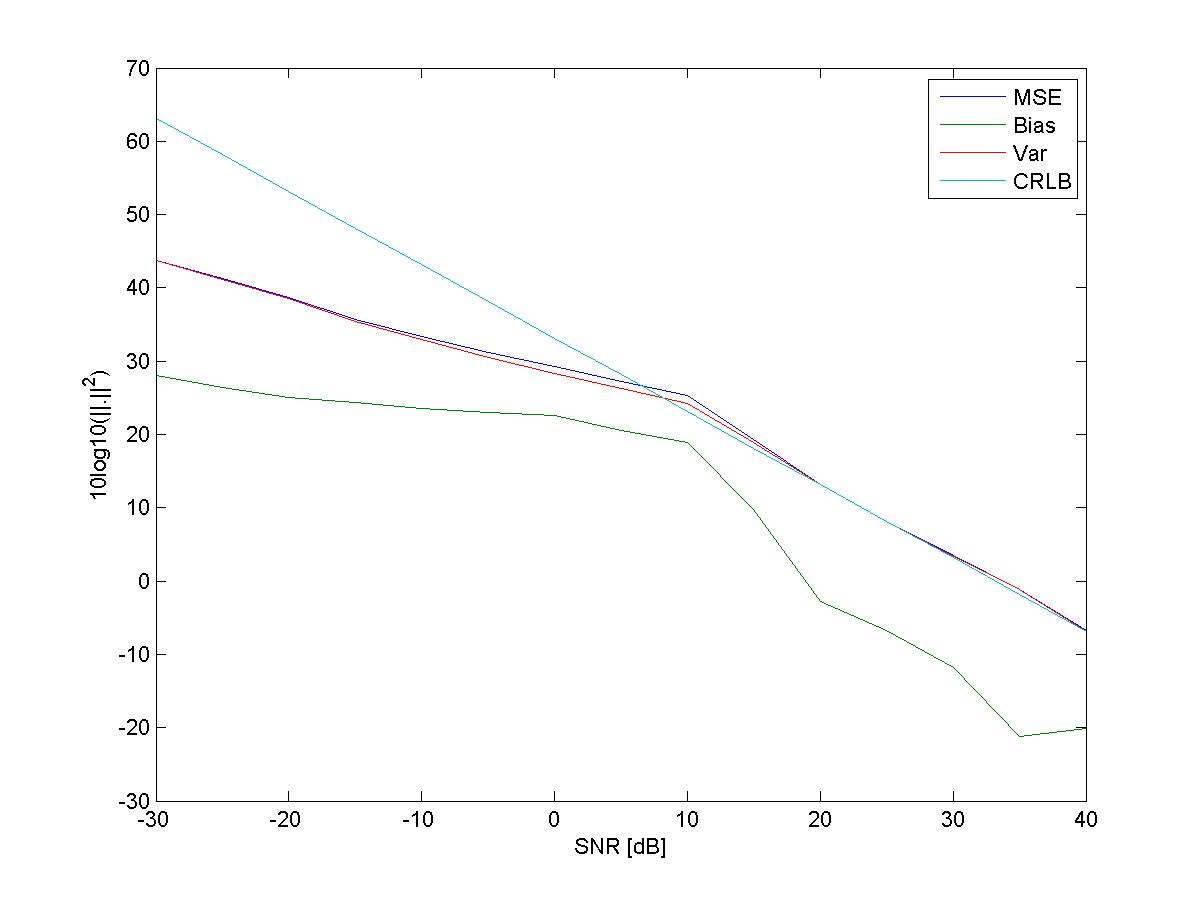}
\caption{Bias and variance components of the mean-square error and CRLB bounds for $m=800$ (top plot), $m=600$ (middle plot),
 and $m=400$ (bottom plot) when $n=100$. \label{fig2}}
\end{figure}

In Figure \ref{fig2} we plot the bias and variance components of the mean-square error for the same results
in Figure \ref{fig1}. Note the bias is relatively small. The bulk of mean-square error is due to estimation variance.

The mean number of iterations varied between 35 and 50, with a higher value for lower $SNR$ (-30dB) and a smaller value for
higher $SNR$ (40dB). 

\begin{figure}[htb]
\includegraphics[width=130mm,height=130mm]{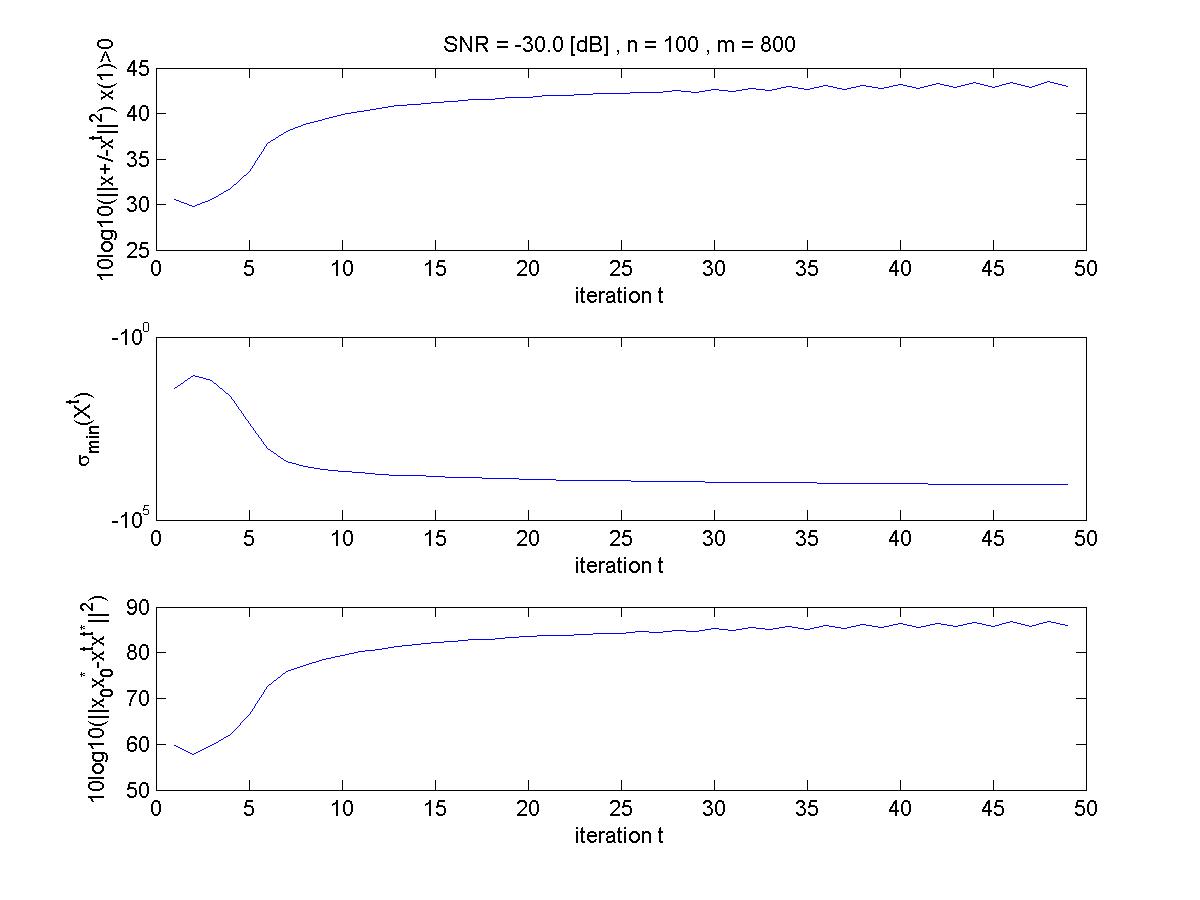}
\includegraphics[width=130mm,height=60mm]{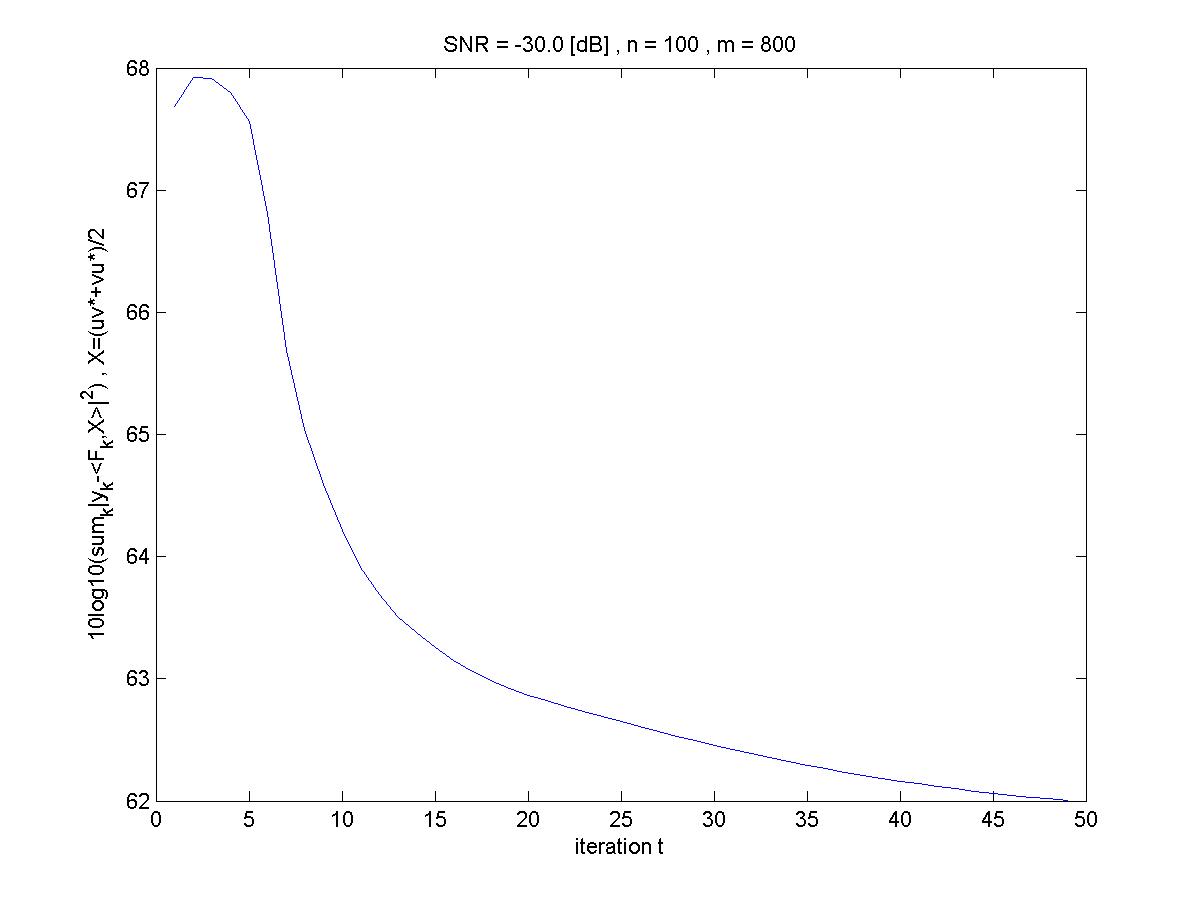}
\caption{Traces for $m=800$, $n=100$ and $SNR=-30dB$: \label{fig3}}
\end{figure}

\begin{figure}[htb]
\includegraphics[width=130mm,height=130mm]{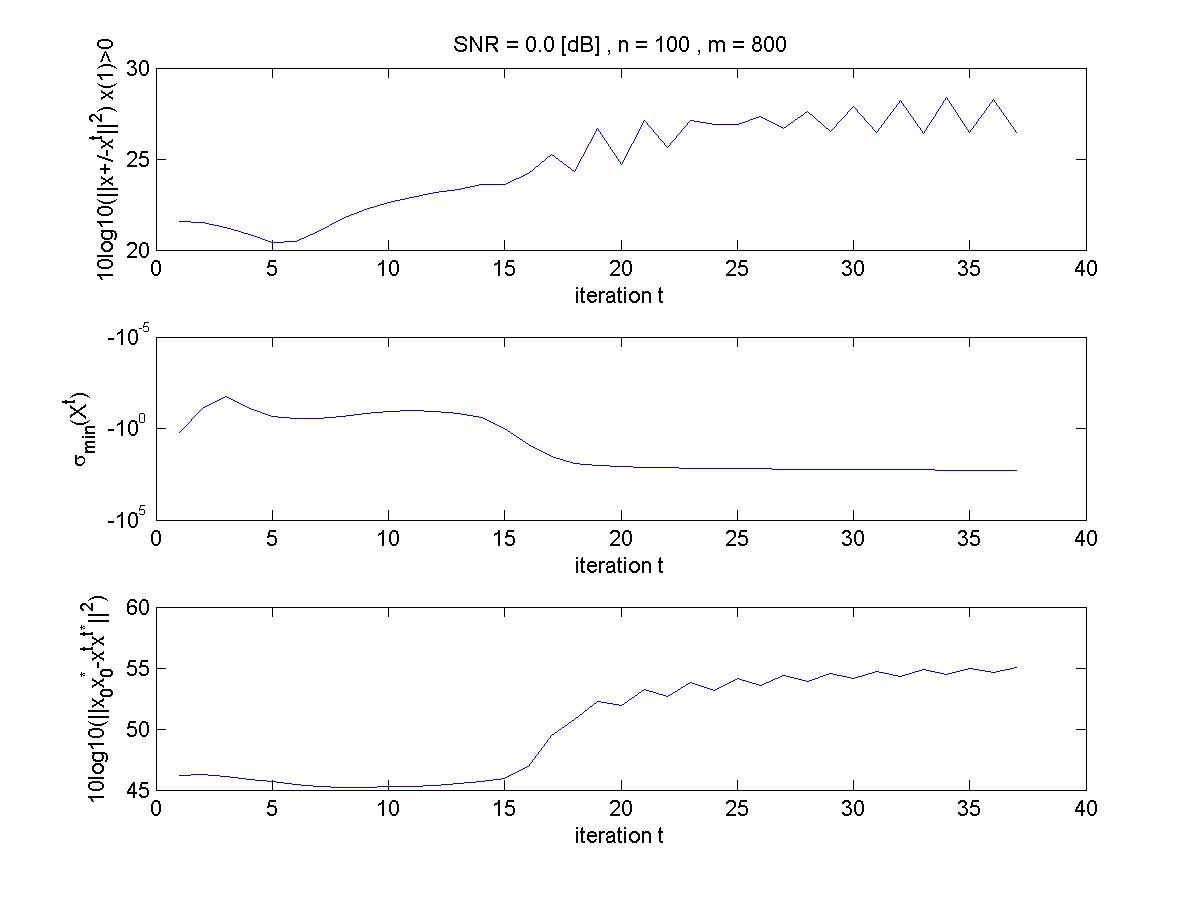}
\includegraphics[width=130mm,height=60mm]{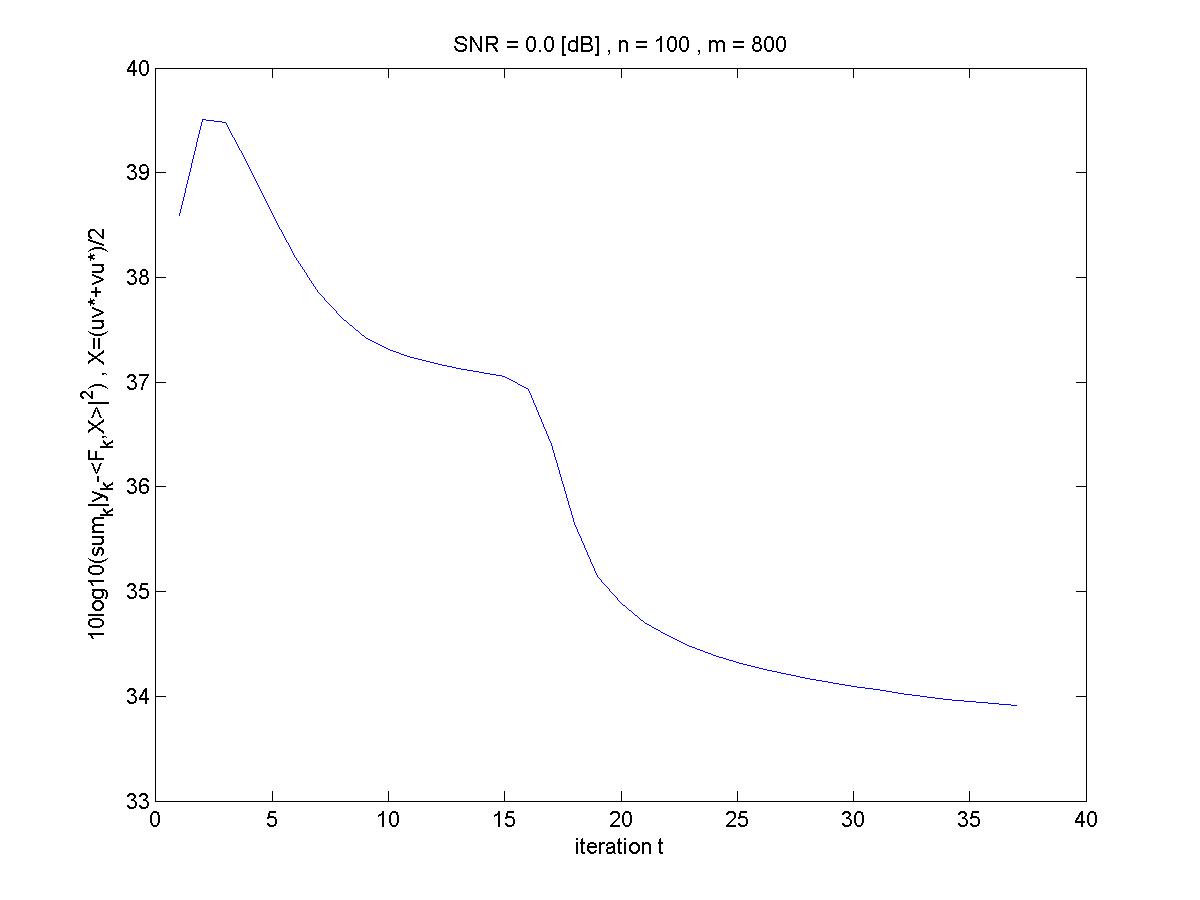}
\caption{Traces for $m=800$, $n=100$ and $SNR=0dB$: \label{fig4}}
\end{figure}

\begin{figure}[htb]
\includegraphics[width=130mm,height=130mm]{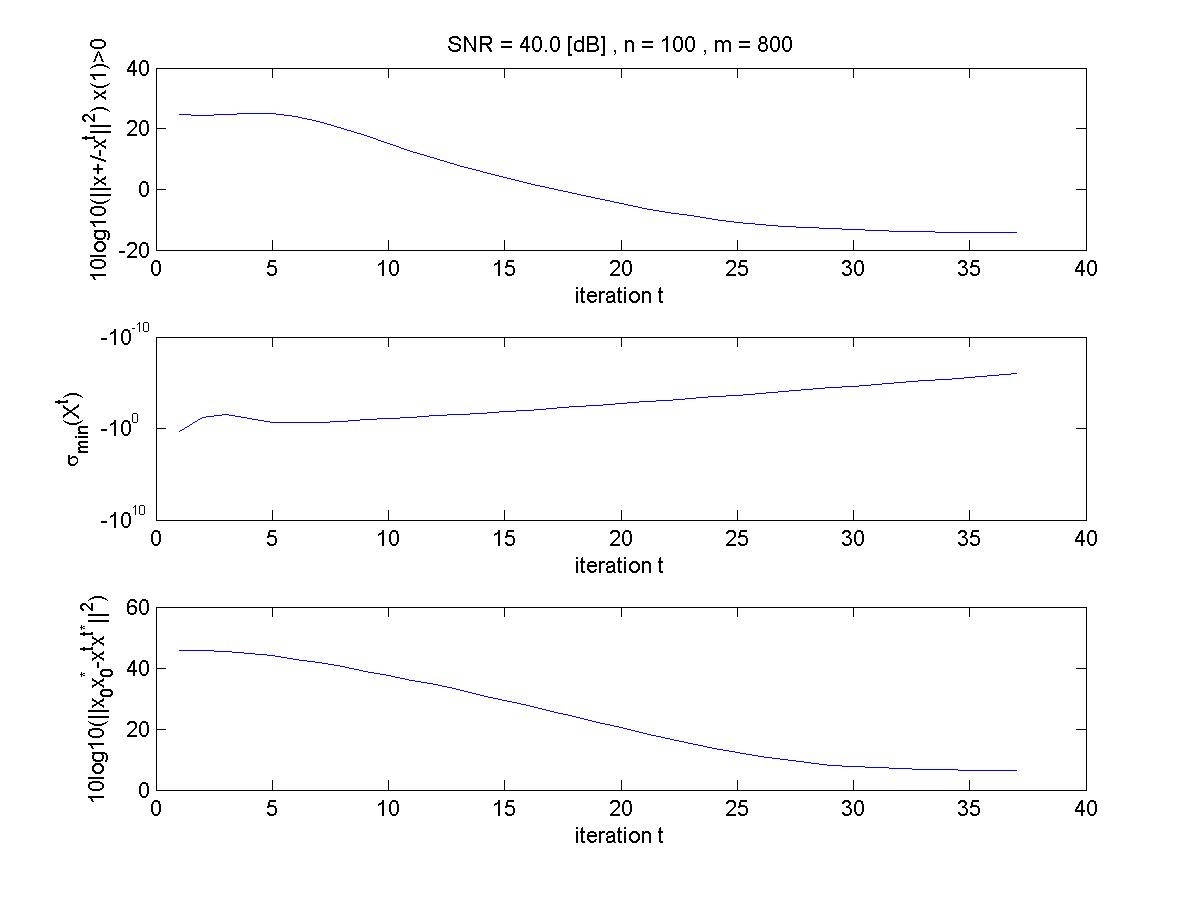}
\includegraphics[width=130mm,height=60mm]{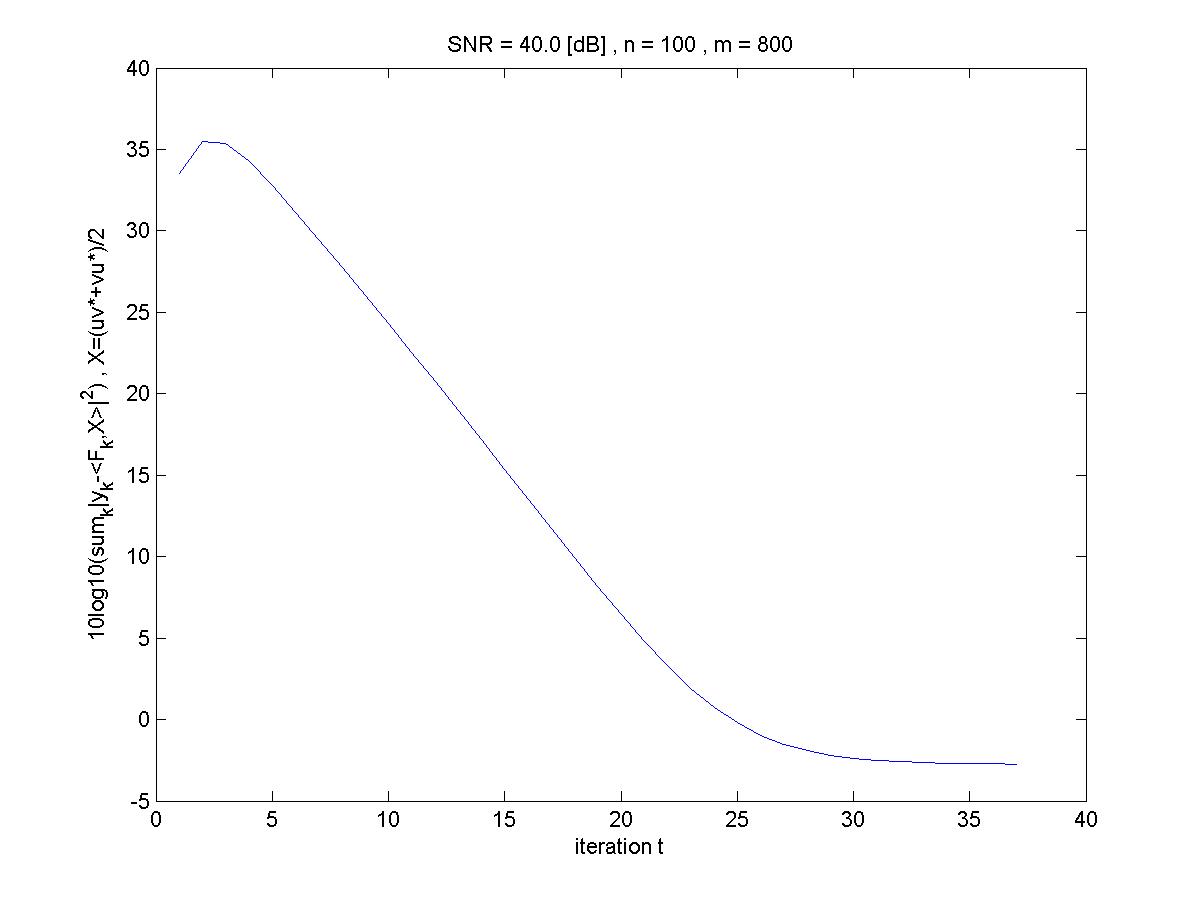}
\caption{Traces for $m=800$, $n=100$ and $SNR=40dB$: \label{fig5}}
\end{figure}

For $m=800$ and three values of $SNR$ (-30dB, 0dB, and 40dB) Figures \ref{fig3}-\ref{fig5}  plot traces of a particular noise realization (different for each SNR). 
In each figure, we plot four characteristics: top graph plots the estimation error of the clean signal $10\log_{10}\norm{x^{(t)}-x}^2$;
second graph presents the smallest eigenvalue of $X^{(t)}=\outp{x^{(t)}}{x^{(t-1)}}$ (which is negative); third graph contains the current estimation error of the clean
rank-one, that is $10\log_{10}\norm{X^{(t)}-xx^*}_2^2$ (Frobenius norm); bottom graph plots $J_0=10\log_{10}\norm{y-\Ac(X^{(t)})}^2$ .

\begin{figure}[htb]
\includegraphics[width=130mm,height=70mm]{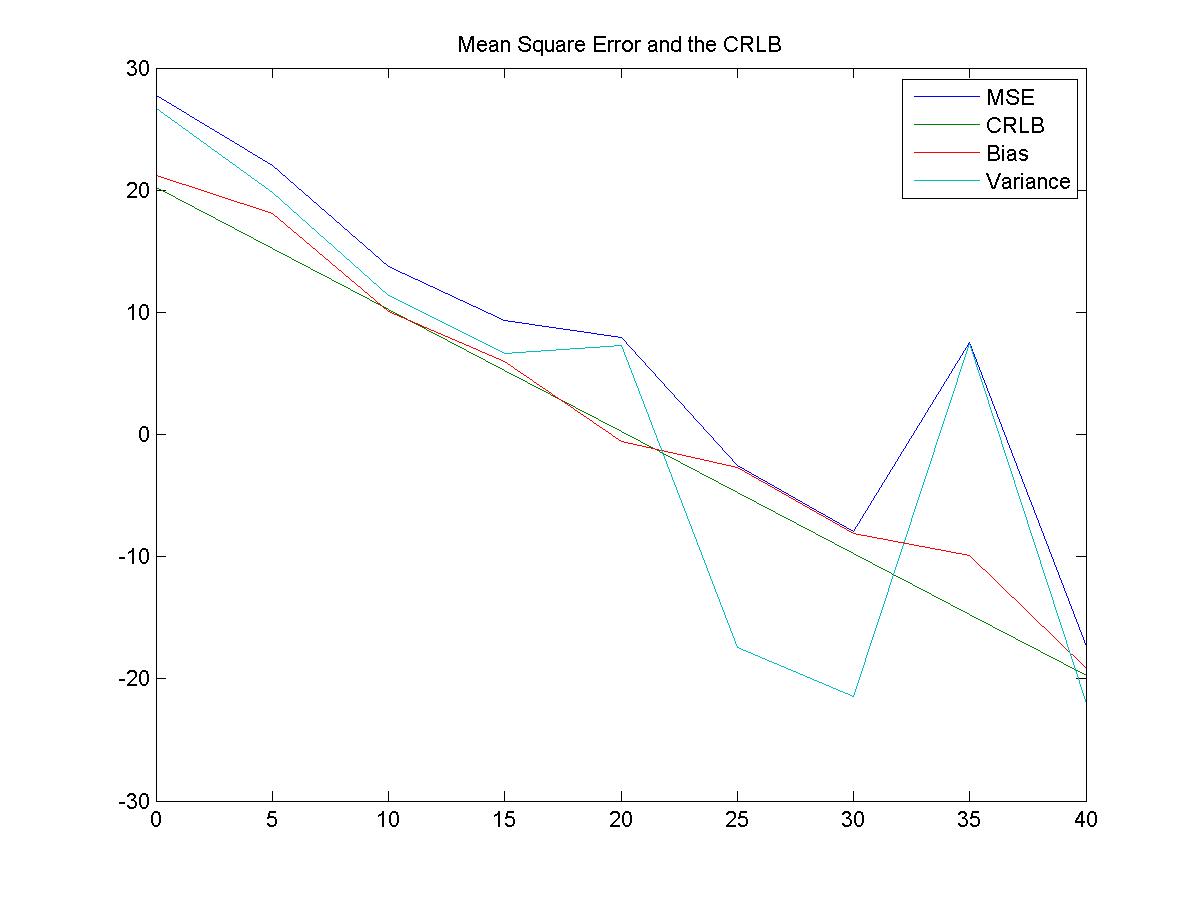}
\includegraphics[width=130mm,height=70mm]{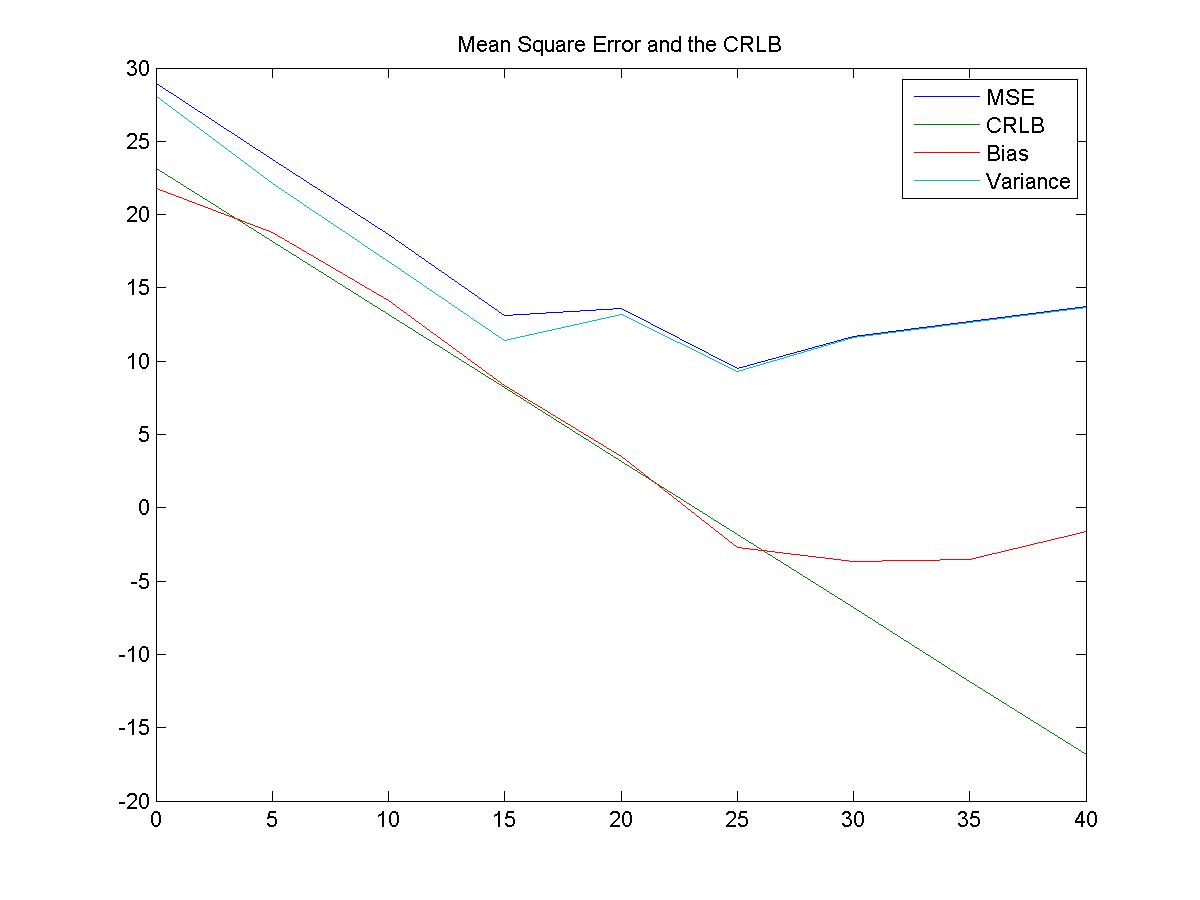}
\includegraphics[width=130mm,height=70mm]{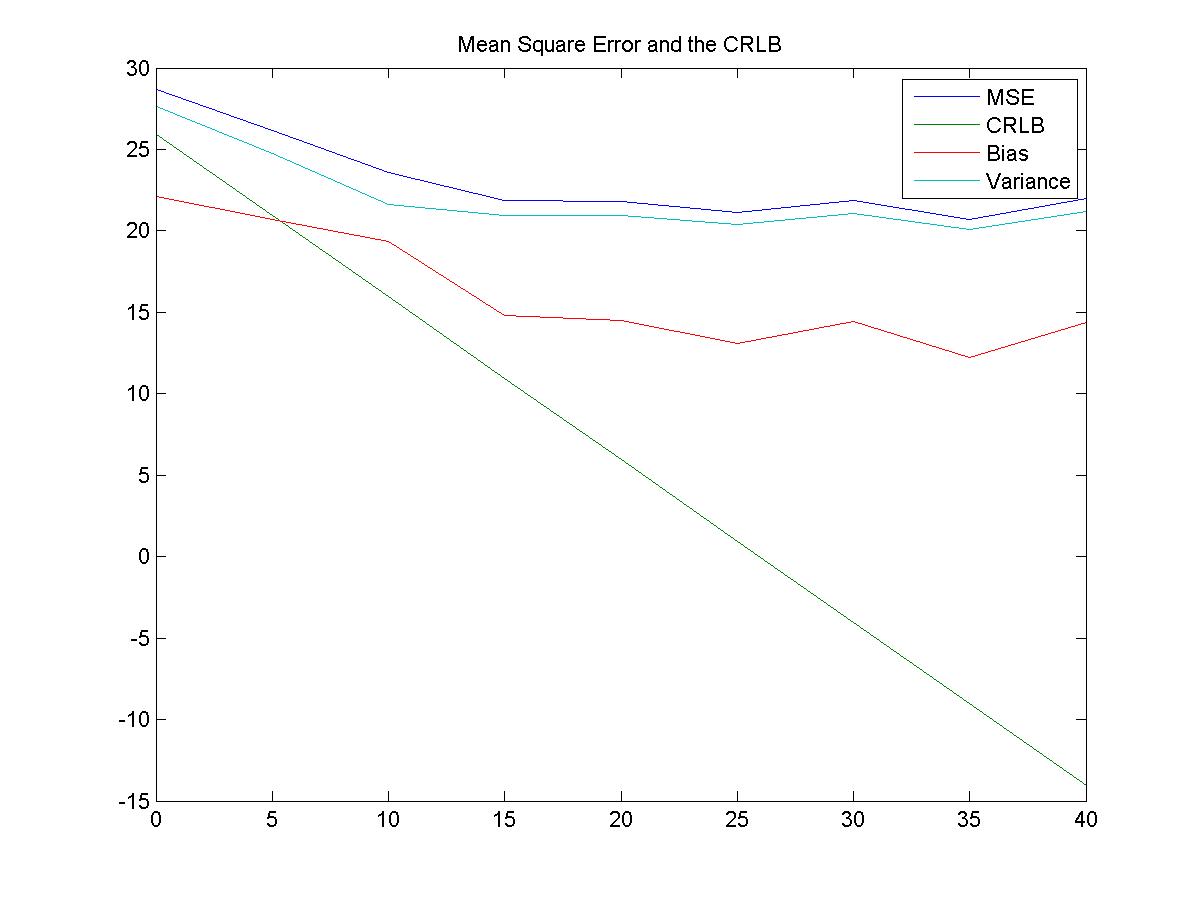}
\caption{Mean-Square Error and CRLB bounds for $m=800$ (top plot), $m=600$ (middle plot),
 and $m=400$ (bottom plot) when $n=100$ and random initialization of $x^0$. \label{fig6}}
\end{figure}

We also analyzed the algorithm sensitivity to initialization. Instead of using the eigenpair (\ref{eq:beta0},\ref{eq:mu0}) we initialize by a random vector $x^0$
together with the eigenvalue (\ref{eq:mu0}). Results for the three values of redundancy are presented in Figure \ref{fig6} for 0dB to 40dB range of $SNR$.

\section{Conclusions\label{sec7}}
Novel necessary conditions for complex valued signal reconstruction from magnitudes of frame coefficients have been presented.
Deterministic stability bounds (Lipschitz constants) and stochastic performance bounds (Cramer-Rao lower bound) have been presented. 
The entire analysis has been done canonically, that is independent of a particular choice of basis.
Then an optimization algorithm based on the least-square error has been proposed and analyzed.
The algorithm performance has been compared to the theoretical lower bound given by the Cramer-Rao inequality. 
Remarkably the algorithm performs very well on a large range of SNR. In particular, for high SNR, it seems to converge to the
correct signal every time. This behavior suggests the algorithm presented here is able to
track the global minimum of (\ref{eq:JJ}) very well. A future study shall analyze this tracking hypothesis.

\section*{Acknowledgments}

The author was partially supported by NSF under DMS-1109498 and DMS-1413249 grants. The author thanks the Erwin Schr\"{o}dinger Institute for the
hospitality shown during the special workshop on ''Phase Retrieval'' in October 2012. Some of the results obtained here were presented at that
workshop and later at the Workshop on ''Phaseless Reconstruction'', UMD, February 2013. The author also thanks Bernhard Bodmann, Jameson Cahill, Martin Ehler, 
Boaz Nadler, Oren Raz, and Yang Wang  for fruitful discussions. He also thanks the annonymous referees for their helpful comments and careful reading of the first draft. Additionally he is grateful to Friedrich Philipp \cite{Fritz} for his comments and for pointing out several errors in the first draft, in particular the dimension of $\Soneone$ is Lemma \ref{l3.2.3}. 
Last but not least, the author thanks the anonymous referees for their comments and their patience with an earlier draft of this paper.

\end{document}